\documentclass[10pt,final,journal]{IEEEtran}

\usepackage{amsmath}
\usepackage{amssymb}
\usepackage{multirow}
\usepackage{algorithm, algpseudocode}
\usepackage{graphicx}
\usepackage[caption=false,subrefformat=parens,labelformat=parens,font=footnotesize]{subfig}
\usepackage{color}
\usepackage{url}
\newcommand{\beq}{\begin{equation}}
\newcommand{\eeq}{\end{equation}}
\newcommand{\bseq}{\begin{subequations}}
\newcommand{\eseq}{\end{subequations}}

\newcommand{\be}{\begin{enumerate}}
\newcommand{\ee}{\end{enumerate}}
\newcommand{\bi}{\begin{itemize}}
\newcommand{\ei}{\end{itemize}}

\newcommand{\cD}{{\cal D}}

\newcommand{\cP}{{\cal P}}

\newcommand{\cR}{{\cal R}}
\newcommand{\cS}{{\cal S}}
\newcommand{\cT}{{\cal T}}

\DeclareMathOperator*{\tfor}{\ for\ }

\DeclareMathOperator*{\argmax}{arg\,max}

\newcommand{\ub}[1]{\overline{#1}}
\newcommand{\prob}[1]{\Pr(#1)}

\newcommand{\omax}[1]{\max_{#1} \,}

\newcommand{\oargmax}[1]{\argmax_{#1} \,}

\newcommand{\norm}[1]{\left\lVert{#1}\right\rVert}
\newcommand{\ip}[2]{\left\langle{#1},{#2}\right\rangle}
\newcommand{\set}[1]{\left\{{#1}\right\}}
\newcommand{\setc}[2]{\left\{{#1}\left\lvert\,{#2}\right.\right\}}
\newcommand{\bm}[1]{\begin{bmatrix}#1\end{bmatrix}}

\def\<{\left<}
\def\>{\right>}
\def\({\left(}
\def\){\right)}
\def\mr{\multirow}
\def\mc{\multicolumn}

\def\SpaRSA{SpaRSA}
\def\MATPOWER{{MATPOWER}}
\def\MATLAB{{\sc Matlab}}

\def\sjwresolved#1{}

\def\noprint#1{}

\newcommand{\tcell}[2][c]{\begin{tabular}[#1]{@{}c@{}}#2\end{tabular}}

\begin{document}

\title{PMU Placement for Line Outage Identification via Multiclass Logistic
  Regression}

\author{Taedong~Kim and Stephen~J.~Wright%
\thanks{This work was supported by a DOE grant
  subcontracted through Argonne National Laboratory Award 3F-30222, and
  National Science Foundation Grant DMS-1216318.}%
\thanks{T. Kim and S.~J. Wright are with the
  Computer Sciences Department, 1210 W. Dayton Street,
  University of Wisconsin, Madison, WI 53706, USA
  (e-mails: tdkim@cs.wisc.edu and swright@cs.wisc.edu).}} %

\maketitle

\begin{abstract}
We consider the problem of identifying a single line outage in a power
grid by using data from phasor measurement units (PMUs). When a line
outage occurs, the voltage phasor of each bus node changes in response
to the change in network topology. Each individual line outage has a
consistent ``signature,'' and a multiclass logistic regression (MLR)
classifier can be trained to distinguish between these signatures
reliably. We consider first the ideal case in which PMUs are attached
to every bus, but phasor data alone is used to detect outage
signatures. We then describe techniques for placing PMUs selectively
on a subset of buses, with the subset being chosen to allow
discrimination between as many outage events as possible. We also
discuss extensions of the MLR technique that incorporate explicit
information about identification of outages by PMUs measuring line
current flow in or out of a bus. Experimental results with
synthetic 24-hour demand profile data generated for 14, 30, 57 and
118-bus systems are presented.
\end{abstract}


\section{Introduction}\label{sec:intro}
In recent years, phasor measurement units (PMUs) have been introduced
as a way to monitor power system networks. Unlike the more
conventional Supervisory Control and Data Acquisitions (SCADA) system,
whose measurements include active and reactive power and voltage
magnitude, PMUs can provide accurate, high-sampling-rate, synchronized
measurements of voltage phasor. There has been much ongoing research
on how the real-time measurement information gathered from PMUs can be
exploited in many areas of power system studies, including system
control and state estimation.

In this paper, we study the use of PMU data in detecting topological
network changes caused by single-line outages, and propose techniques
for determining optimal placement of a limited number of PMU devices
in a grid, so as to maximize the capability for detecting such
outages. Our PMU placement approach can also be used as a tie-breaker
for the other types of strategies that have multiple optimal
solutions, for example, maximum observability problems.

Knowledge of topological changes as a result of line failure can be
critical in deciding how to respond to a blackout. Rapid detection of
such changes can enable actions to be taken that reduce risks of
cascading failures that lead to large-scale blackouts. One of the main
causes of the catastrophic Northeast blackout of 2003 was faulty
topological knowledge of the grid following the initial failures (see
\cite{BlackOut2003}).

Numerous approaches have been proposed for identifying line outages
using PMU device measurements. In \cite{TatO08, TatO09}, phasor angle
changes are measured and compared with expected phasor angle
variations for all single- or double-line outage scenarios. Support
vector machines (SVM) were proposed for identification of single-line
outages in \cite{AbdME12}. A compressed-sensing approach was applied
to DC power balance equations to find sparse topological changes in
\cite{ZhuG12}, while a cross-entropy optimization technique was
considered in \cite{CheLM14}.
Since the approaches in \cite{ZhuG12} or \cite{CheLM14} use the
linearized DC power flow models to represent a power system, their
line outage identification strategies rely only on changes to phase
angles. Voltage magnitude measurements from PMUs, which also provide
useful information for monitoring a power system, are ignored.  Our
use of the AC power flow model allow both more accurate modeling of
the system and more complete exploitation of the available data.

The key feature that makes line outage identification possible is that
voltage phasor measurements reported by PMUs are different for
different line-outage scenarios. Our approach aims to distinguish
between these different ``signatures'' by using a multiclass logistic
regression (MLR) model. The model can be trained by a convex
optimization approach, using standard techniques. The coefficients
learned during training can be applied during grid operation to detect
outage scenarios. Our approach could in principle be applied to
multiline outages too, but since the problem dimension is much larger
in such cases --- the number of possible outage scenarios is much
greater --- it is no longer practical. Second, even when trained only
to recognize single-line outages, our classifier is useful in
multiline outage situations on large grids, when the coupling between
the lines is weak (as discussed in \cite[Section~2.2]{ZhuS14}). In
other words, many multiline outage cases can be decomposed into
single-line outage events on different parts of the grid.

Because of the expense of installation and maintenance, PMUs can
reasonably be installed on just a subset of buses in a grid. We
therefore need to formulate an {\em optimal placement problem} to
determine the choice of PMU locations that gives the best information
about system state.  Several different criteria have been proposed to
measure quality of a given choice of PMU locations. One of the most
popular criteria is
to place PMUs to maximize the number of nodes in the system that can
be observed directly, either by a PMU located at that node or an
adjacent connected node \cite{ChaKE09}. Another criterion is
quality of state estimation results.  In this approach, one can use
PMU measurements alone, or combine them with traditional SCADA
measurement to decide the optimal PMU deployment (see for example
\cite{KekGW12}). Other criteria and techniques for locating PMUs
optimally are discussed in the review papers \cite{ManKG11, YuiEC11}.

For the case in which line outage identification is used as a
criterion for PMU placement, we mention \cite{ZhaGP12,ZhuG12,ZhuS14}.
In \cite{ZhaGP12}, the authors use pre-computed phase angles as outage
signatures and attempt to find the optimal PMU locations by
identifying a projection (by setting to zero the entries which are not
selected as PMU locations) that maximizes the minimum distance
in $p$-norm of the projected signatures.  The problem is formulated as
an integer program (IP) and a greedy heuristic and branch-and-bound
approach are proposed.  PMU placement for the line outage
identification method discussed in \cite{ZhuG12} is studied in
\cite{ZhuS14}.  A non-convex mixed-integer nonlinear program (MINLP)
is formulated, leading to a linear programming convex relaxation.
Again, a greedy heuristic and a branch-and-bound algorithm is
suggested as a solution methodology.

We build our PMU placement formulation on our MLR model for
single-line outage detection, by adding nonsmooth ``Group LASSO''
regularizers to the MLR objective and applying several heuristics.

The rest of this paper is organized as follows: In
Section~\ref{sec:line.outage}, the line outage identification problem
is described along with the multiclass logistic regression (MLR)
formulation. The problem of PMU placement to identify a line outage is
described in Section~\ref{sec:PMU.placement}, and we describe the
group-sparse heuristic and its greedy variant used to formulate and
solve this problem. Experimental results on synthetically generated
data are presented in Section~\ref{sec:result}. A conclusion appears
in Section~\ref{sec:conclusion}.

In supplementary material, we describe an extension that makes use of
explicit line outage information. This model uses the fact that when a
PMU is attached to a line, it can detect by direct observation when
that line fails, and has no need to rely on the indirect evidence of
voltage changes. As expected, performance improves when such
additional information is used, though as we show in this paper, very
good results can still be obtained even when it is ignored.

\section{Line Outage Identification}\label{sec:line.outage}

We describe an approach that uses changes in voltage phasors measurements at
PMUs to detect single-line outages in the power grid. As in
\cite{TatO08,ZhuG12}, we assume that the fast dynamics of the system are well
damped and voltage measurements reflect the quasi-static equilibrium that is
reached after the disruption. We use a quasi-steady state AC power flow model
(see e.g. \cite[Chapter~10]{BerV99}) as a mapping from time varying load
variation (and line outage events) to the polar coordinate ``outputs'' of
voltage magnitude and angle.

PMUs report phasor measurements with high frequency, and changes in
voltage due to topology changes of the power grid tend to be larger
than the variation of voltage phasor during normal operation (for
example, demand fluctuation that occurs during the sampling time
period). We construct signature vectors from these voltage changes
under the various single-line outage situations, and use them to train
a classifier.

\begin{figure*}[t]\centering\footnotesize
 \raisebox{90pt}{\rotatebox{90}{Voltage Angle (rad)}}
 \includegraphics[width=.95\textwidth]{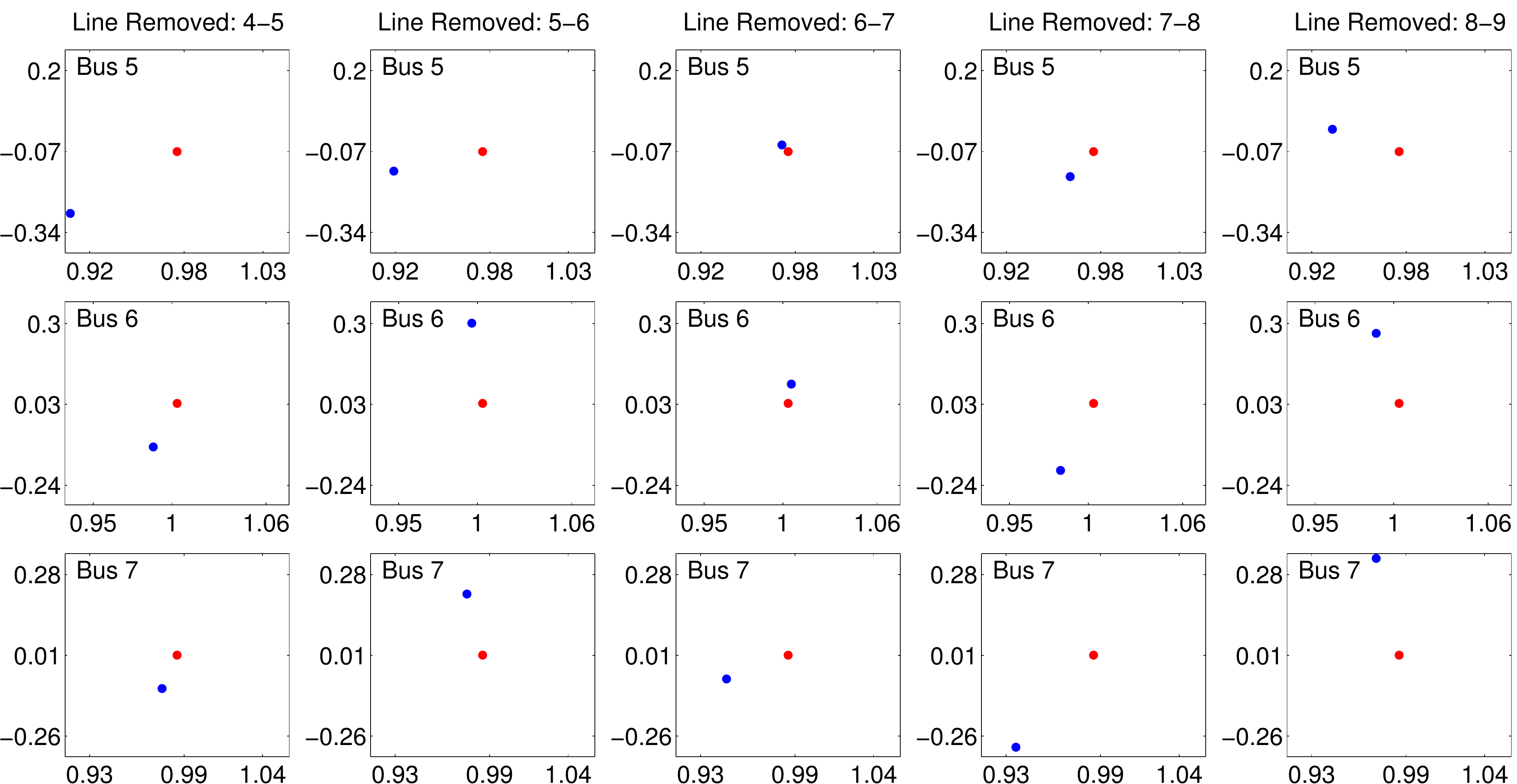}
 \mbox{\hspace{26pt} Voltage Magnitude (p.u.)}
  \caption{9-Bus System: Voltage changes caused by single line outage on buses
  5, 6 and 7. (Voltage magnitude (p.u) - Phase angle (rad))}
  \label{fig:9bus.phasor.change}
\vspace*{-3pt}
\end{figure*}

Figure~\ref{fig:9bus.phasor.change} shows an example of voltage
changes for a 9-bus system ({\tt case9.m} from
\MATPOWER{}~\cite{ZimMT11}) on different line outage scenarios. The
failure of lines connecting buses 4-5, 5-6, 6-7, 7-8, and 8-9 is
considered as possible scenarios (columns in the figure) whose voltage
phasors at buses 5, 6, and 7 (rows in the figure) observed.  In each
plot, $x$-axis shows voltage magnitude and $y$-axis shows voltage
phase angle at the bus. The red dots in each plot indicate the voltage
phasor when there is no line failure and the blue dots are the voltage
values under the specified failure scenario. We observe that voltage
values at these buses change in distinctive ways under different line
outage scenarios. It therefore seems realistic to expect that by
comparing voltage phasor data, gathered before and after a failure
event, we can identify the failed line reliably. We now describe the
multinomial logistic regression model for determining the outage
scenario.

\subsection{Multinomial Logistic Regression Model}

Multinomial logistic regression (MLR) is a machine-learning approach
for multiclass classification. In our application, examples of voltage
phasor changes under each outage scenario are used to train the
classifier by determining parameter values in a set of parametrized
functions. Once the parameters have been found, these functions
determine the likelihood of a given set of phasor changes as being
indicative of each possible failure scenario.

Suppose that there are $K$ possible outcomes (classes) labelled as $i
\in \set{1,2,\dotsc,K}$ for a given vector of observations $X$. In the
multinomial logistic regression model, the probability of a given
observation $X$ has an outcome $Y$ (one of the $K$ possibilities
$i\in\set{1,2,\cdots,K}$) is given by the following formula:
\begin{align}\label{eq:prob.Y}
  \prob{Y=i|X} &:= \frac{e^{\ip{\beta_i}{X}}}
                      {\sum_{k=1}^Ke^{\ip{\beta_k}{X}}}
  \tfor i=1,2,\cdots,K,
\end{align}
where $\beta_1, \beta_2, \dotsc, \beta_K$ are regression coefficients,
whose values are obtained during the training process. Note that there
is one regression coefficient $\beta_i$ for each outcome $i \in
\{1,2,\dotsc,K\}$. Once values of the coefficients $\beta_i$, $i \in
\{1,2,\dotsc,K\}$ have been determined, we can predict the outcome
associated with a given feature vector $X$ by evaluating
\[
  k^* = \oargmax{k\in\set{1,2,\dotsc,K}}\prob{Y=k|X},
\]
or equivalently,
\beq\label{eq:max.outcome}
  k^* = \oargmax{k\in\set{1,2,\dotsc,K}} \ip{\beta_k}{X}.
\eeq

Training of the regression coefficients $\beta_1, \cdots, \beta_K$ can
be performed by maximum likelihood estimation. The training data
consists of $M$ pairs $(X_1,Y_1), (Y_2, Y_2), \dotsc, (X_M, Y_M)$,
each consisting of a feature vector and its corresponding
outcome. Given formula \eqref{eq:prob.Y}, the a posteriori likelihood
of observing $Y_1,Y_2,\dotsc,Y_M$ given the events
$X_1,X_2,\dotsc,X_M$ is
\beq\label{eq:likelihood} \prod_{i=1}^M P(Y=Y_i|X_i)
\,=\, \prod_{i=1}^M \(\frac{e^{\ip{\beta_{Y_i}}{X_i}}}
      {\sum_{k=1}^Ke^{\ip{\beta_k}{X_i}}}\).
\eeq
By taking $\log$ of \eqref{eq:likelihood}, we have log-likelihood
function
\[
  f(\beta) :=
    \sum_{i=1}^M\(\ip{\beta_{Y_i}}{X_i} -
    \log\sum_{k=1}^Ke^{\ip{\beta_k}{X_i}}\),
\]
where the matrix $\beta$ is obtained by arranging the coefficient
vectors as $\bm{\beta_1 & \beta_2 & \dotsc & \beta_K}$. The maximum
likelihood estimate $\beta^*$ of regression coefficients is obtained
by solving the following optimization problem:
\beq\label{eq:MLE}
\beta^* = \oargmax{\beta} f(\beta).
\eeq
This is a smooth convex problem that can be solved by fairly standard
techniques for smooth nonlinear optimization, such as L-BFGS
\cite{LiuN89}. Note that $f(\beta) \le 0$ for all $\beta$.

If the training data is separable, the value of $f(\beta)$ can be made
to approach zero arbitrarily closely by multiplying $\beta$ by an
increasingly large positive value (see \cite{KriCF05}). To recover
meaningful values of $\beta$ in this case, we can solve instead the
following regularized form of \eqref{eq:MLE}:
\beq\label{eq:MLE.p}
\beta^* = \oargmax{\beta} f(\beta) - \tau w(\beta)
\eeq
where $\tau>0$ is a penalty parameter and $w(\beta)$ is a (convex)
penalty function of the coefficient $\beta$. The penalized form can
also be used to promote some kind of structure in the solution
$\beta^*$, such as sparsity or group-sparsity.  This property is key
to our PMU placement formulation, and we discuss it further in
Section~\ref{sec:PMU.placement}.

Training of the MLR model, via solution of \eqref{eq:MLE} or
\eqref{eq:MLE.p}, can be done offline, as described in the next
subsection. Once the model is trained (that is, the coefficients
$\beta_i$, $i=1,2,\dotsc,K$ have been calculated), classification can
be done via \eqref{eq:max.outcome}, at the cost of multiplying the
matrix $\beta$ by the observed feature vector $X$, an operation that
can be done in real time.

\subsection{Training Data: Observation Vectors and Outcomes}\label{sec:data.format}
In our MLR model for line outage identification problems, the
observation vector $X_j$ is constructed from the change of voltage
phasor at each bus, under a particular outage scenario. The
corresponding outcome is the index of the failed line.

Suppose that a power system consists of $N$ buses, all equipped with
PMUs that report the voltage values periodically. Let $(V_i,
\theta_i)$ and $(V_i', \theta_i')$, $i=1,2,\cdots, N$, be two phasor
measurements obtained from PMU devices, one taken before a possible
failure scenario and one afterward. The observation vector $X$ which
describes the voltage phasor difference is defined to be
\beq \label{eq:X}
  X = \bm{\Delta V_1 & \cdots & \Delta V_N &
        \Delta\theta_1 & \cdots & \Delta\theta_N}^T
\eeq
where $\Delta V_i=V_i'-V_i$ and $\Delta\theta_i = \theta_i'-\theta_i$,
for $i=1,2,\dotsc,N$. If we assume that the measurement interval is
small enough that loads and demands on the grid do not change
significantly between measurements, we would expect the entries of $X$
to be small, unless an outage scenario (leading to a topological
change to the grid) occurred. Some such outages would lead to failure
of the grid. More often, feasible operation can continue, but with
significant changes in the voltage phasors, indicated by large
components of $X$.

The training data $(X_j,Y_j)$ can be assembled by a considering a
variety of realistic demand scenarios for the grid, solving the AC
power flow equations for each possible outage scenario (setting the
value of $Y_j$ according to the index of that failure), then setting
$X_j$ to be the shift in voltage phasor that corresponds to that
scenario. The phasor shifts for a particular scenario change somewhat
as the pattern of loads and generations changes, so it is important to
train the model using a sample of phasor changes under different
realistic patterns of supply and demand.

The observation vector can be extended to include additional
information beyond the voltage phasor information from the PMUs, if
such information can be gathered easily and exploited to improve the
performance of the MLR approach. For example, the system operator may
be able to monitor the power generation level $G$ (expressed as a
fraction of the long-term average generation) that is injected to the
system at the same time points at which the voltage phasor
measurements are reported.  If included in the observation vector,
this quantity might need to scaled so that it does not dominate the
phasor difference information. Also, a constant entry can be added to
the observation vector to provide more flexibility for the
regression. The extended observation vector thus has the form
\beq \label{eq:Xbar}
  \ub{X} = \bm{\Delta V_1 & \cdots & \Delta V_N &
            \Delta\theta_1 & \cdots & \Delta\theta_N & \rho G & \rho}^T
\eeq
where $\rho$ is a scaling factor that approximately balances the
magnitudes of all entries in the vector. (Note that since $G$ is not
too far from $1$, it is appropriate to use the same scaling factor for
the last two terms.)  The numerical experiments in
Section~\ref{sec:result} make use of this extended observation vector.

\section{PMU Placement}\label{sec:PMU.placement}

As we mentioned in Section \ref{sec:intro}, installing of PMUs at all
buses is impractical. Indeed, if it were possible to do so,
single-line outage detection would become a trivial problem, as each
outage could be observed directly by PMU measurements of line current
flows in or out of a bus; there would be no need to use the
``indirect'' evidence provided by voltage phasor changes. In this
section we address the problem of placing a limited number of PMUs
around the grid, with the locations chosen in a fashion that maximizes
the system's ability to detect single-line outages.  This PMU
placement problem selects a subset of {\em buses} for PMU placement,
and assumes that PMUs are placed to monitor voltage phasors at the
selected buses.

A naive approach is simply to declare a ``budget'' of the number of
buses at which PMU placement can take place, and consider all possible
choices that satisfy this budget. This approach is of course
computationally intractable except for very small cases.  Other
possible approaches include a mixed-integer nonlinear programming
formulation \cite{ZhaGP12,ZhuS14}, but this is very hard to solve in
general. In this paper, we use a regularizer function $w(\beta)$ in
\eqref{eq:MLE.p} to promote the a particular kind of sparsity
structure in the coefficient matrix $\beta$. Specifically, A group
$\ell_1$-regularizer is used to impose a common sparsity pattern on
all columns in the coefficient matrix $\beta$, with nonzeros occurring
only in locations corresponding to the voltage magnitude and phase
angle changes for a subset of buses. The numerical results show that
approaches based on this regularizer give reasonable performance on
the PMU placement problem.

\subsection{Group-Sparse Heuristic (GroupLASSO)}\label{sec:gsh}
Let $\cP$ be the set of indices in the vector $X$, that is
$\cP=\set{1,2,\cdots,|X|}$. Consider $S$ mutually disjoint subsets of
$\cP$, denotes $\cP_1, \cP_2, \dotsc, \cP_S$. For each $s \in \cS :=
\set{1,2,\cdots,S}$, define $q_s(\beta)$ as follows:
\[
  q_s(\beta) = \norm{[\beta]_{\cP_s}}_{F}
             = \sqrt{\sum_{i\in\cP_s}\sum_{k=1}^K (\beta_{ik})^2}
\]
where $[\beta]_{\cP_s}$ is the submatrix of $\beta$ constructed by
choosing the rows whose indices are in $\cP_s$, $\norm{\cdot}_F$ is
the Frobenius norm, and $\beta_{ik}$ is the $(i,k)$ entry of matrix
$\beta$ (thus, $\beta_{ik}$ is the $i$th entry of the coefficient
vector $\beta_k$ in \eqref{eq:max.outcome} and \eqref{eq:likelihood})
The value of $q_s(\beta)$ is the $\ell_2$-norm over the entries of
matrix $\beta$ which are involved in group $s$. For our observation
vectors $X$ \eqref{eq:X} and $\ub{X}$ \eqref{eq:Xbar}, we can choose
the number of groups $|\cS|$ equal to the number of buses $N$, and set
\beq \label{eq:Ps}
  \cP_s = \{ s, s+N \}, \quad s=1,2,\dotsc,N.
\eeq
Thus, if bus $s$ is ``selected'' in the placement problem, the
coefficients associated with $\Delta V_s$ and $\Delta \theta_s$ are
allowed to be nonzeros. Buses that are not selected need not of course
be instrumented with PMUs, because the coefficients in $\beta$ that
correspond to these buses are all zero. Note that for the
extended vector $\ub{X}$, we do not place the last two entries (the
constant and the total generation) into any group, as we assume that
these are always ``selected'' for use in the classification process.

For any subset $\cR$ of $\cS$, we define a group-$\ell_1$-regularizer
$w_{\cR}(\beta)$  to be the sum of $q_s(\beta)$ for $s\in\cR$,
that is,
\[
  w_{\cR}(\beta) = \sum_{s\in\cR} q_s(\beta).
\]
Setting $\cR=\cS$, the penalized form \eqref{eq:MLE.p} with
$w=w_{\cS}$ can be used to identify a group-sparse solution:
\beq\label{eq:MLE.group.l1}
  \omax{\beta} f(\beta) - \tau w_{\cS}(\beta).
\eeq
With an appropriate choice of the parameter $\tau$, the solution
$\beta^*$ of \eqref{eq:MLE.group.l1} will be group-row-sparse, that
is, the set $\setc{s\in\cS}{q_s(\beta^*)\neq 0}$ will have
significantly fewer than $|\cS|$ elements.  Given a solution $\beta^*$
of \eqref{eq:MLE.group.l1} for some value of $\tau$, we could define
the $r$-sparse solution as follows (for a given value of $r$, and
assuming that the solution of \eqref{eq:MLE.group.l1} has at least $r$
nonzero values of $q_s(\beta^*)$):
\beq\label{eq:MLE.l1.S}
 \cR^* := \oargmax{\cR:|\cR|=r,\cR\subset\cS} w_{\cR}(\beta^*).
\eeq
Since the minimizer $\beta^*$ of \eqref{eq:MLE.group.l1} is biased due
to the presence of the penalty term, we should not use the submatrix
extracted from $\beta^*$ according to the selected group $\cR^*$ as
the regression coefficients for purposes of multiclass
classification. Rather, we should solve a reduced, unpenalized version
of the problem in which just the coefficients from sets $\cP_s$ that
were {\em not} selected are fixed at zero. That is, we define a {\em
  debiased} solution $\tilde{\beta}^*$ corresponding to $\cR^*$ as follows:
\begin{align}\label{eq:MLE.debias}
  \omax{\beta}f(\beta) \;\;
& \mbox{subject to $\beta_{ik}=0$ for all $(i,k)$ with } \\
\nonumber
& \mbox{ $k=1,2,\dotsc,K$ and  $i\in \cP_s$ for some $s \notin \cR^*$.}
\end{align}

\subsection{Greedy Heuristic}

The regularization approach can be combined with a greedy strategy, in
which groups are selected one at a time, with each selection made by
solving a regularized problem.
Suppose that $\cR^{l-1}$ is set of selected groups after $l-1$
iterations of the selection heuristic. The problem solved at iteration
$t$ of the heuristic to choose the next group is
\beq\label{eq:MLE.l1.greedy}
\hat{\beta}^l =  \oargmax{\beta} f(\beta) - \tau w_{\cS \setminus \cR^{l-1}}(\beta).
\eeq
The next group $s^l$ is obtained from $\hat{\beta}^l$ as follows:
\[
  s^l = \oargmax{s\in \cS \setminus \cR^{l-1}} q_s(\hat{\beta}^l),
\]
and we set $\cR^l=\cR^{l-1}\cup\set{s^l}$.  Note that we do {\em not}
penalize groups in $\cR^{l-1}$ that have been selected already, in
deciding on the next group $s^l$.
After choosing $r$ groups by this process, the debiasing step is
performed to find the best maximum likelihood estimate for the sparse
observation. Algorithm~\ref{alg:PMU.placement.greedy} describes this
greedy approach. Note that the initial set of groups $\cR^0$ might not
be empty since we can use additional information that is independent
from the PMU measurement, if available.

\begin{algorithm}\small
\caption{Greedy Heuristic}\label{alg:PMU.placement.greedy}
\begin{algorithmic}[1]
\Require
  \Statex Choose an initial set of groups: $\cR^0$.
  \Statex Parameter $\tau$, $r$.
\Ensure
  \Statex $\cR^r$: Set of groups after selecting $r$ groups.
  \Statex $\tilde{\beta}^r$: Maximum likelihood estimate for $r$-group
  observation.
\medskip
\For {$l=1,2,\cdots,r$}
  \State Solve \eqref{eq:MLE.l1.greedy} with $\cR^{l-1}$ for $\hat{\beta}^l$.
  \State $s^l\gets\oargmax{s\in\cS \setminus \cR^{l-1}} q_s(\beta^l)$
  \State $\cR^l\gets \cR^{l-1}\cup\set{s^l}$
\EndFor
\State Solve \eqref{eq:MLE.debias} with $\cR^r$ to get $\tilde{\beta}^r$.
$\Comment{debiasing}$
\end{algorithmic}
\end{algorithm}

The major advantage of this approach is that redundant observations
are suppressed by already-selected, non-penalized observations at each
iteration.  We will give more details in discussing the experimental
results in Section~\ref{sec:result}.

\section{Numerical Experiments}\label{sec:result}

Here we present experimental results for the approaches proposed
above.  The test sets considered here are based on the power system
test cases from \MATPOWER{} (originally from \cite{PSTCA}), with
demands altered to generate training and test sets for the MLR
approach.

\subsection{Synthetic Data Generation}
\begin{figure*}\centering
\subfloat[Base Case Demand (Considered as Average)]
  {\includegraphics[width=.3\textwidth]{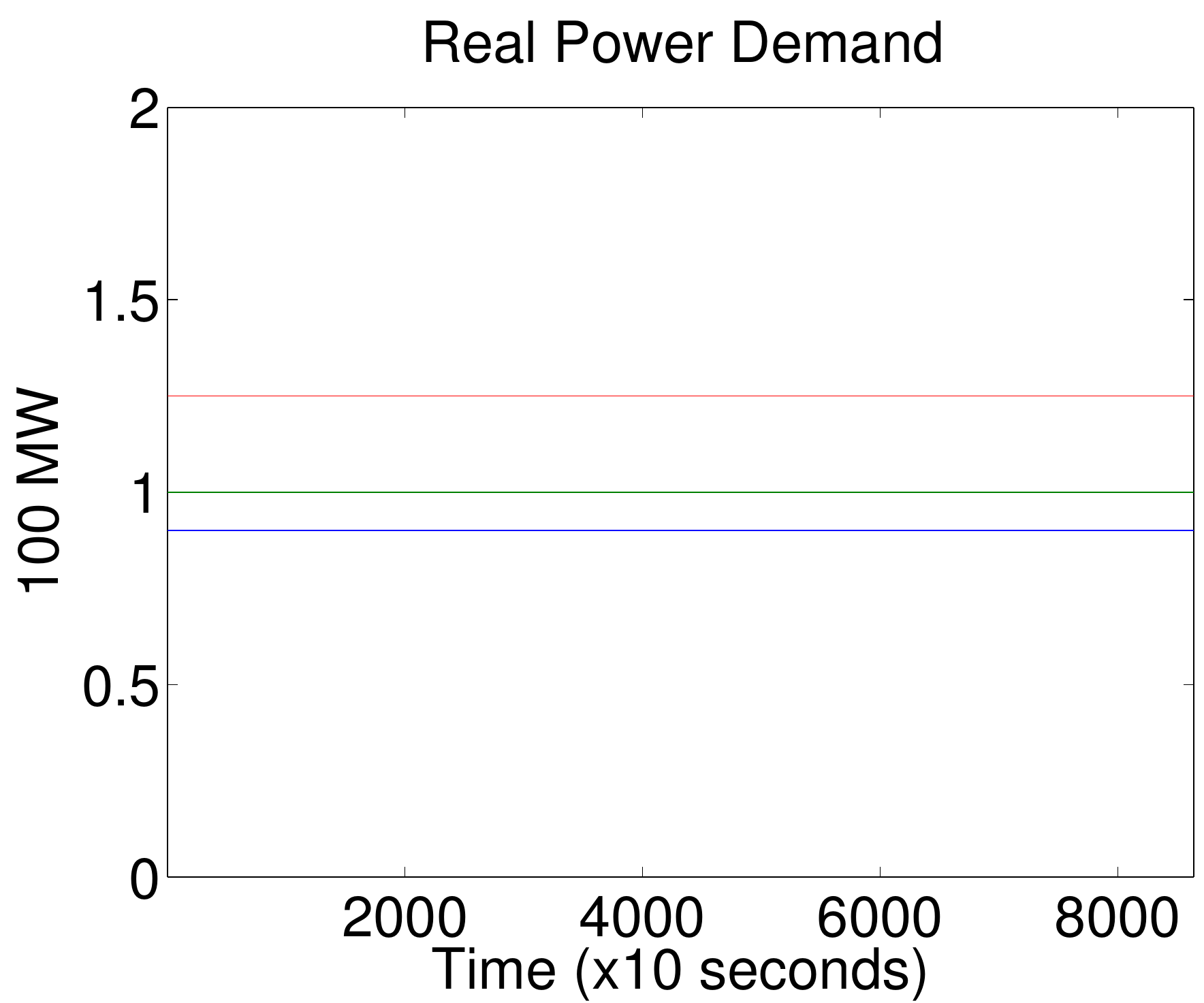}
    \label{fig:syn.data.avg}}\qquad
\subfloat[The Ratio Generated by Stochastic Process]
  {\includegraphics[width=.3\textwidth]{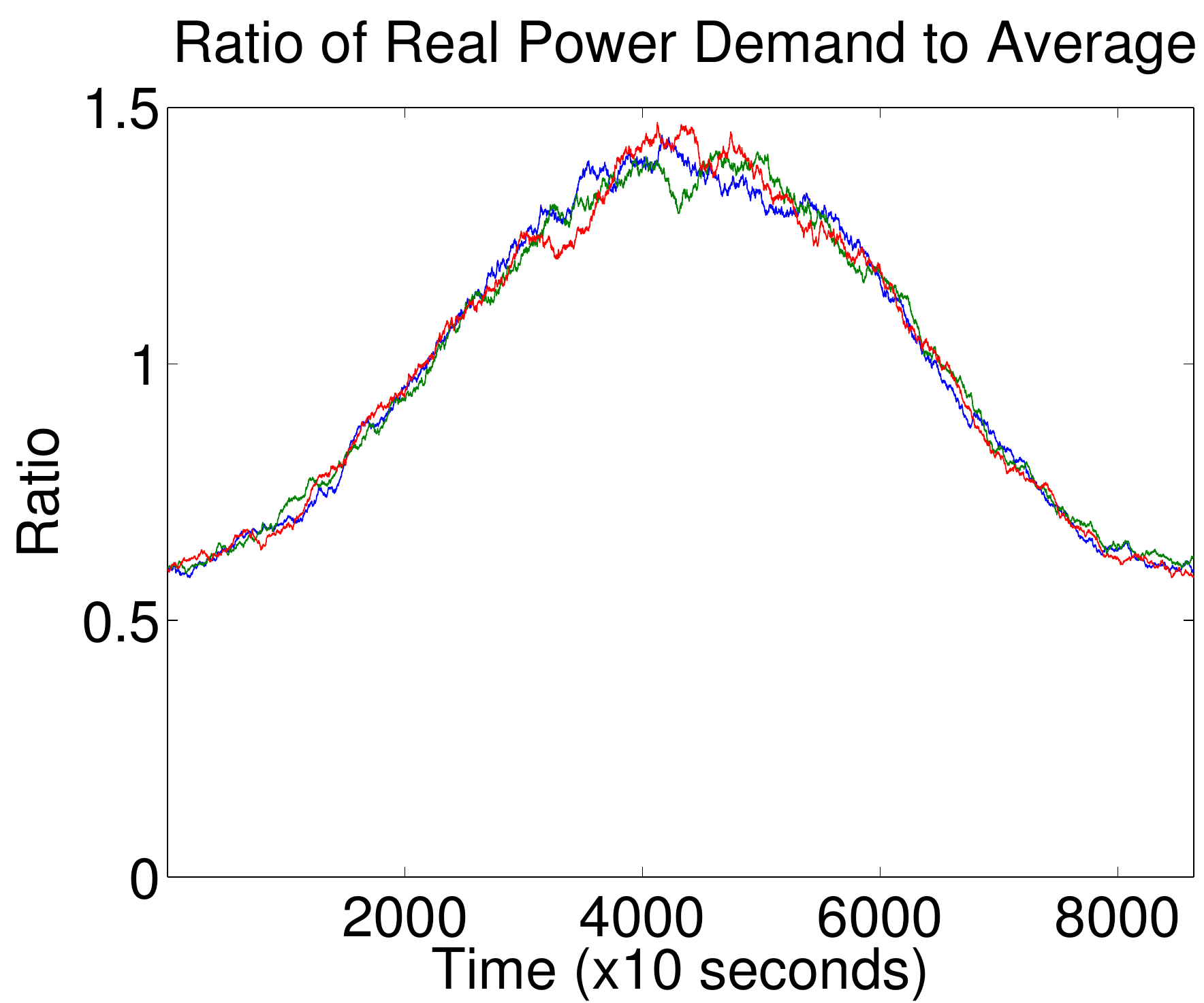}
    \label{fig:syn.data.ratio}}\qquad
\subfloat[Generated 24-Hour Demand Profile]
  {\includegraphics[width=.3\textwidth]{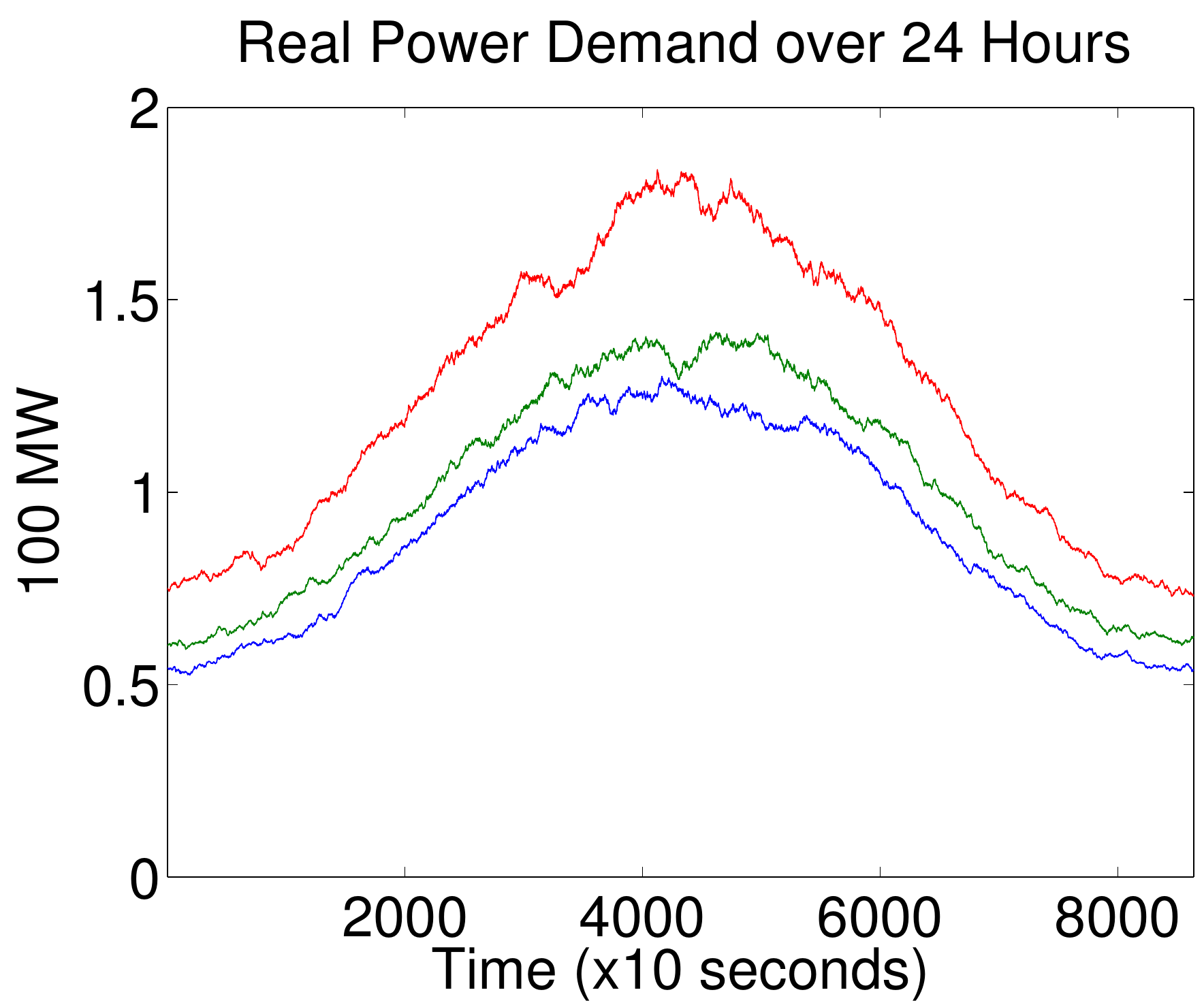}
    \label{fig:syn.data.demand}}
\caption{Generating Synthetic Demand Data by A Stochastic Process}\label{fig:syn.data}
\vspace*{-3pt}
\end{figure*}

Since the data provided from IEEE test case archive \cite{PSTCA} is a
single snapshot of the states of power systems, we extend them to a
synthetic 24-hour demand data cycle by using a stochastic process, as
follows.
\be
  \item Take the demand values given by the IEEE test case archive as
    the average load demand over 24-hours.
  \item Generate the demand variation profile by using an additive
    Ornstein-Uhlenbeck process as described in \cite{PerKA11},
    separately and independently on each demand bus.
  \item Combine the average demand and the variation ratio to obtain
    the 24-hour load demand profile for the system.
\ee
Figure~\ref{fig:syn.data} shows demand data generated by this procedure at
three demand buses in the 9-bus system ({\tt case9.m}) from \MATPOWER{}.
Figure \subref*{fig:syn.data.avg} shows the data drawn from the \MATPOWER{}
file, now taken to be a 24-hour average.  Figure~\subref*{fig:syn.data.ratio}
shows the ratio generated by the additive Ornstein-Uhlenbeck process, and
Figure \subref*{fig:syn.data.demand} shows the products of the average and
ratio.
Since the power injected to the system needs to increase
proportionally to the total demands, all power generation is
multiplied by the average of the demand ratios.  This average of
ratios is used as the generation level $G$ for the observation vector
$\ub{X}$ defined by \eqref{eq:Xbar}.  The data assumes a 10-second
interval between the measurements, so the total number of time points
in the generated data is $24\times60\times6=8640$.

Once the 24-hour load demand profile is obtained, the AC-power
equations are solved using \MATPOWER{} to calculate the voltage phasor
values at each time point. These phasor values are taken to be the PMU
measurements for a normal operation cycle over a 24-hour period.
\MATPOWER{}'s AC power flow equations solver is also used to evaluate
voltage phasors for each single-line outage scenario that does not
lead to an infeasible system.  (During this process, if there exist
duplicated lines that connect the same pair of buses, they are
considered as a single line, that is, we do not allow only a fraction
of multiple lines that connect the same set of buses to be failed.)
Simulation of single-line failures to generate training data is
necessary because there are typically few instances of actual outages
available for study.  The voltage variation for each line outage at
time $t$ is calculated by subtracting these normal-operation voltages
at timepoint $t-1$ from line outage voltages at time point $t$. (The
10-second interval between measurements is usually sufficient time to
allow transient fluctuations in phasor values to settle down; see
\cite{TatO08}.) This process leads to a number of labeled data pairs
$(X,Y)$ (or $(\ub{X},Y)$) which we can use to train or tune the MLR
classifier.

\begin{table}\centering
\caption{Test Cases from \MATPOWER{}}\label{tbl:test.cases}
\begin{tabular}{|c|c|c|c|c|c|}\hline
\mr{2}{*}{System}  & \MATPOWER{}
  & \mc{2}{c|}{\# of Lines} & Train & Test \\
\cline{3-4}
  size & case & Feas. & Infeas./Dup. & (5) & (50) \\\hline\hline
  14-Bus  &        case14  &   18  &    2  &    90  &   900 \\
  30-Bus  &  case\_ieee30  &   37  &    4  &   185  &  1850 \\
  57-Bus  &        case57  &   67  &   13  &   335  &  3350 \\
 118-Bus  &       case118  &  170  &   16  &   850  &  8500 \\\hline
\end{tabular}
\end{table}

Table~\ref{tbl:test.cases} provides the basic information on the power
systems used for the experiments. The number of lines that are
feasible is given in the column ``Feas.,'', while the number of lines
that are duplicated or that lead to an infeasible power flow problem
when removed from the system is shown in the column
``Infeas./Dup.''. For each {\em feasible} line outage, five equally
spaced samples are selected from the first half (that is, the first
12-hour period) of voltage variation data as training instances.
Fifty samples are selected randomly from the second half of voltage
variation data as test instances. The numbers of training and test
instances are shown in the last two columns of the table.

\subsection{PMUs on All Buses}

\begin{table}\centering
\caption{Line Outage Detection Accuracy on Test Set with PMUs on All Buses.}
\vspace{-15pt}
\subfloat[Based on Probability of Correct Answer]{
  \begin{tabular}{|c||c|c|c||c|c|c|}\hline
  \mr{2}{*}{System} & \mc{3}{c||}{Using $X$} & \mc{3}{c|}{Using $\ub{X}$} \\
  \cline{2-7}
   & $\ge0.9$ & $\ge0.7$ & $\ge0.5$ & $\ge0.9$ & $\ge0.7$ & $\ge0.5$
     \\\hline\hline
  14-Bus &  100\% &  100\% &  100\% &  100\% &  100\% &  100\% \\
  30-Bus & 99.7\% & 99.7\% & 99.7\% &  100\% &  100\% &  100\% \\
  57-Bus & 99.5\% & 99.7\% & 99.8\% & 99.5\% & 99.7\% & 99.8\% \\
 118-Bus & 99.5\% & 99.5\% & 99.5\% & 99.8\% & 99.8\% & 99.8\% \\\hline
  \end{tabular}\label{tbl:result.all_bus.prob}
  }\\
  \subfloat[Based on Ranking of Correct Answer]{
  \begin{tabular}{|c||c|c|c||c|c|c|}\hline
  \mr{2}{*}{System} & \mc{3}{c||}{Using $X$} & \mc{3}{c|}{Using $\ub{X}$} \\
  \cline{2-7}
     & $1$ & $\le2$ & $\le3$ & $1$ & $\le2$ & $\le3$ \\\hline\hline
  14-Bus &  100\% &  100\% &  100\% &  100\% &  100\% & 100\% \\
  30-Bus & 99.7\% &  100\% &  100\% &  100\% &  100\% & 100\% \\
  57-Bus & 99.8\% & 99.9\% & 99.9\% & 99.8\% &  100\% & 100\% \\
 118-Bus & 99.5\% & 99.7\% & 99.7\% & 99.8\% & 99.9\% & 100\% \\\hline
  \end{tabular}\label{tbl:result.all_bus.rank}
  }\\\vspace{5px}
\bi
\item ``Probability'' indicates statistics for the probability
  assigned by the MLR classifier to the actual outage
  event.\\
\item ``Ranking'' indicates whether the actual event was ranked in
  the top 1, 2, or 3 of probable outage events by the MLR
  classifier.
\ei\vspace{-10px}
\end{table}

We present results for line outage detection when phasor measurement
data from all buses is used.
The maximum likelihood estimation problem \eqref{eq:likelihood} with
these observation vectors is solved by L-BFGS algorithm \cite{LiuN89},
coded in \MATLAB. We measure performance of the identification
procedure in two ways. The first measure is based on the probability
assigned by the model to the actual line outage.
Table~\subref*{tbl:result.all_bus.prob} shows the accuracy of the
classifiers according to this measure, for both the original phasor
difference vector $X$ \eqref{eq:X} and the extended vector
$\ub{X}$ \eqref{eq:Xbar}.  Each column shows the percentage of testing
samples for which the probability assigned to the correct outage
exceeds $0.9$, $0.7$ and $0.5$, respectively. The result shows that
the performance of line outage identification is very good, even for
the original observation vector $X$.
For both of $X$ and $\ub{X}$, the accuracy of line outage
identification based on probability $\ge 0.5$ is at least $99\%$.

The second measure is obtained by ranking the probabilities assigned
to each line outage on the test datum, and score a positive mark if
the correct outage is one of the top one, two, or three cases in the
ranking. We see in Table~\subref*{tbl:result.all_bus.rank} that the
actual case appears in the top two in almost every case.


\subsection{PMU Placement}\label{sec:result.PMU.placement}

In this subsection we only consider the extended observation vector
$\ub{X}$ defined by \eqref{eq:Xbar}. We assume too that a PMU is
installed on the reference bus, for purposes of maintaining
consistency in phase angle measurement. We describe in some detail the
performance of the proposed algorithm on the IEEE 57 Bus system,
showing that line-outage identification performance when PMUs are
placed judiciously almost matches performance in the
fully-instrumented case. We then summarize our computational
experience on 14, 30, 57, and 118-bus systems.

For our regularization schemes, we used groups $\cP_s$,
$s=1,2,\dotsc,N$, defined as in \eqref{eq:Ps}. The final two entries
in the extended observation vectors (the average-generation and
constant terms) are not included in any group.


\subsubsection{IEEE 57 Bus System}
\begin{figure*}\centering
\subfloat[$\tau=10^{-5}$]{
  \label{fig:PMU.placement.GroupLASSO.1}
  \includegraphics[width=.328\textwidth]{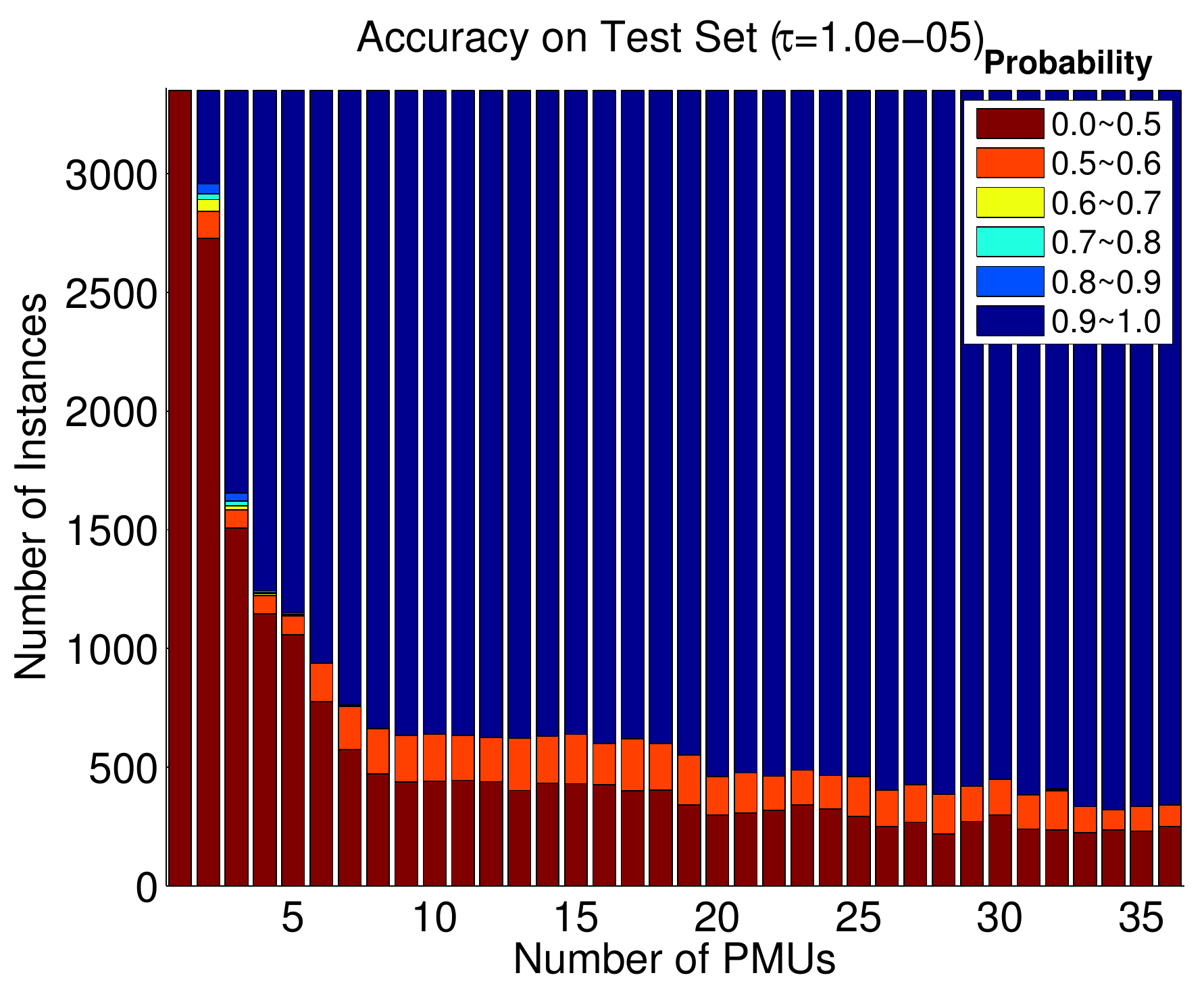}}
\subfloat[$\tau=10^{-3}$]{
  \includegraphics[width=.328\textwidth]{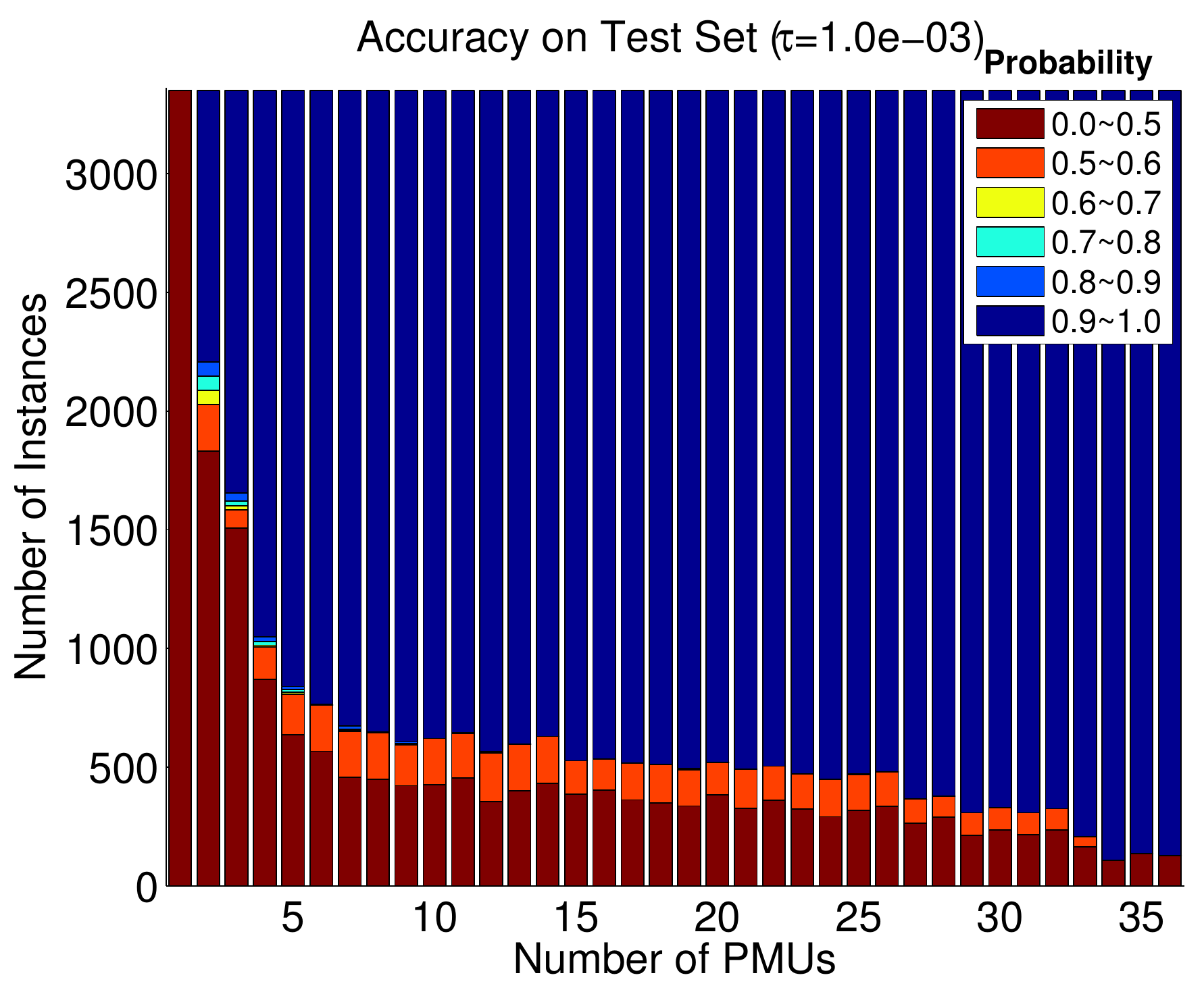}}
\subfloat[$\tau=10^{-1}$]{
  \label{fig:PMU.placement.GroupLASSO.3}
  \includegraphics[width=.328\textwidth]{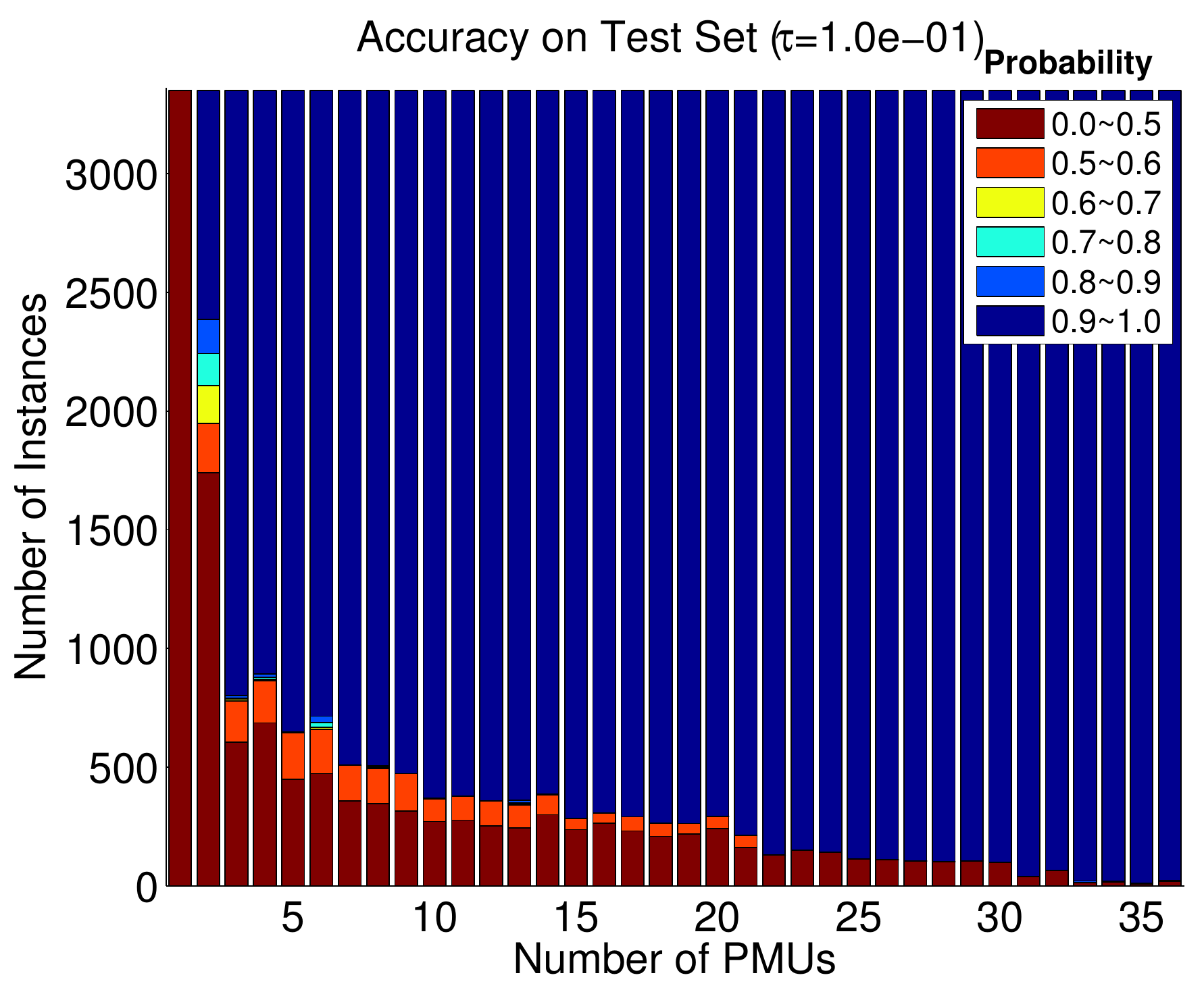}}\\
\caption{Accuracy on Test Set of IEEE 57 Bus System for different values of $\tau$: Group-Sparse
Heuristic}\label{fig:PMU.placement.GroupLASSO}
\subfloat[$\tau=10^{-5}$]{
  \includegraphics[width=.328\textwidth]{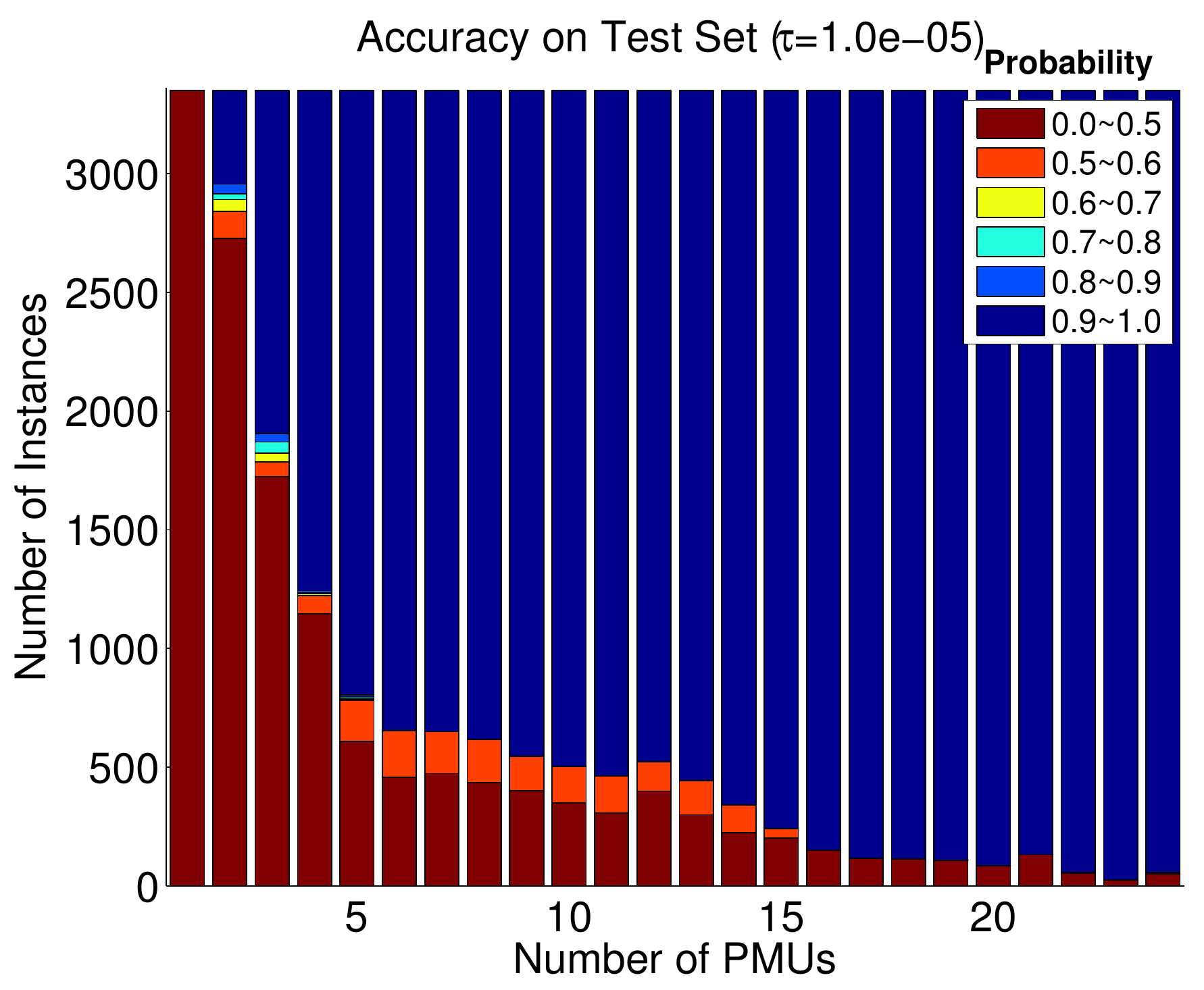}}
\subfloat[$\tau=10^{-3}$]{
  \label{fig:PMU.placement.greedy.2}
  \includegraphics[width=.328\textwidth]{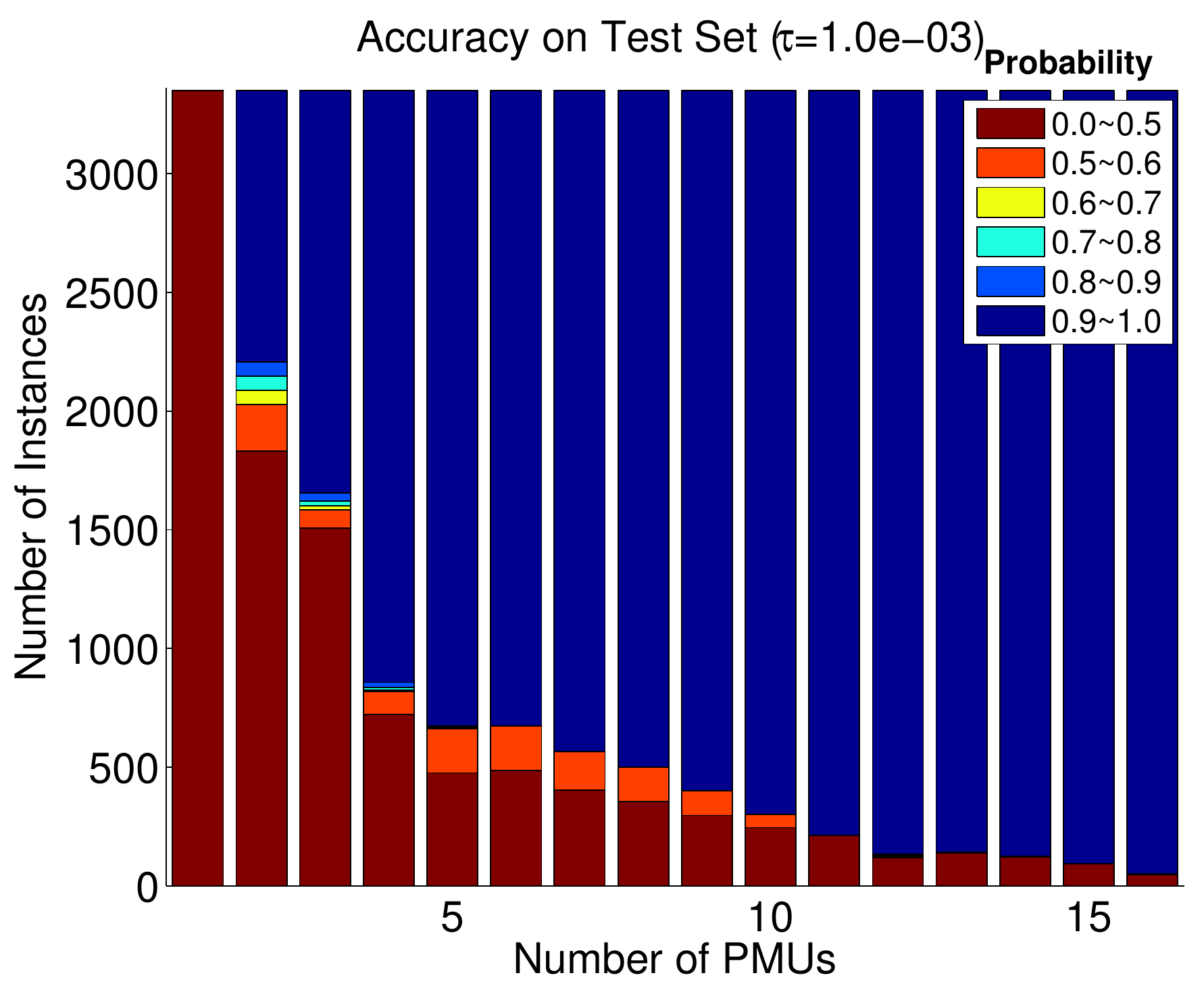}}
\subfloat[$\tau=10^{-1}$]{
  \label{fig:PMU.placement.greedy.3}
  \includegraphics[width=.328\textwidth]{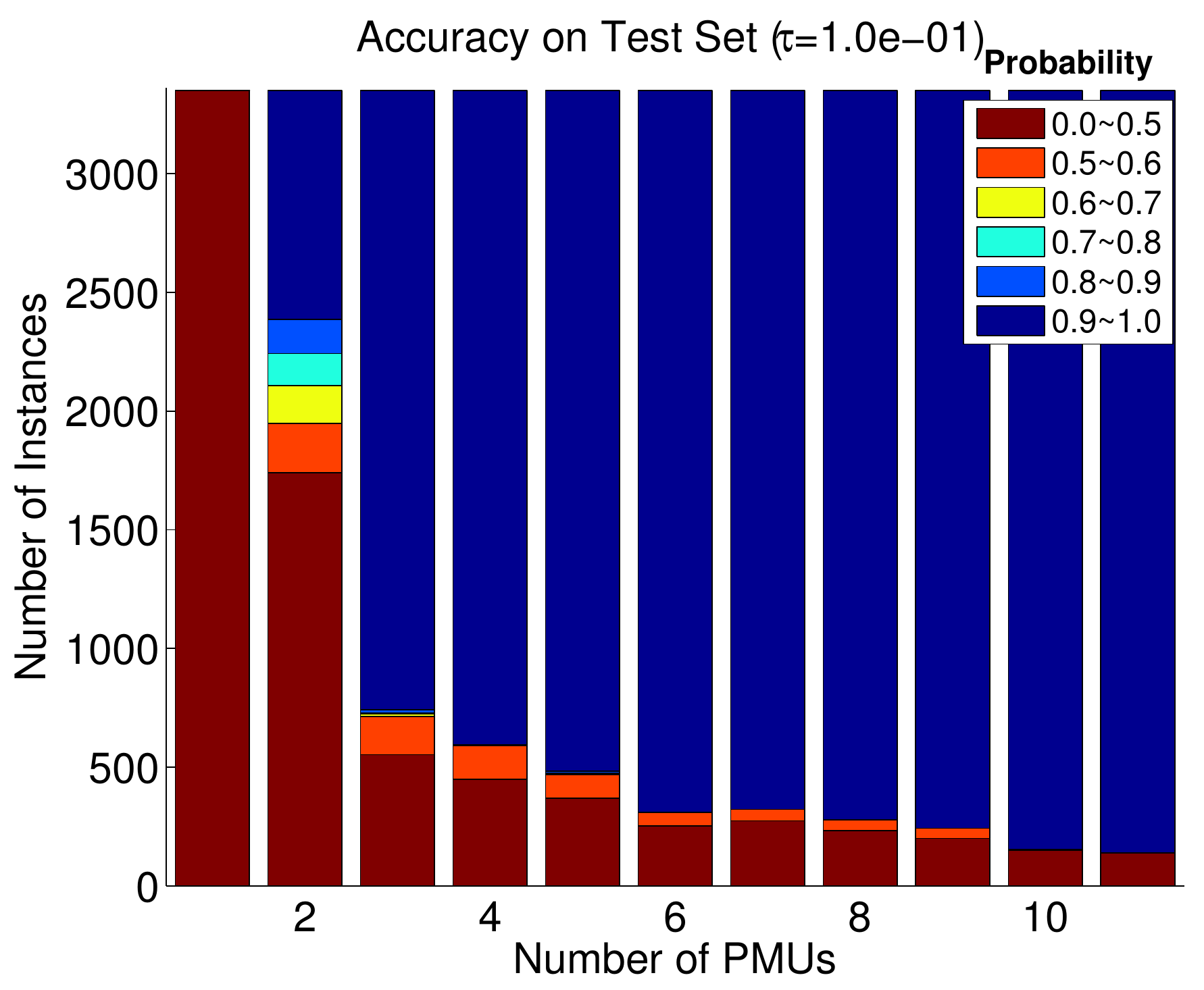}}
\caption{Accuracy on Test Set of IEEE 57 Bus System for different values of $\tau$: Greedy Heuristic}\label{fig:PMU.placement.greedy}
\vspace*{-4pt}
\end{figure*}

We describe here results obtained on the IEEE 57-bus system with two
heuristics discussed in Section~\ref{sec:PMU.placement}: The
GroupLASSO and Greedy Heuristics.

In Figure~\ref{fig:PMU.placement.GroupLASSO}, results for the
GroupLASSO heuristic are displayed for different values of $\tau$. The
$x$-axis indicates the number of PMUs selected by this heuristic. The
$y$-axis indicates the number of test cases for which the true outage
was classified by the heuristic. Each bar is color-coded according to
the probability assigned to the true outage by the MLR
classifier. Blue colors indicate that a high probability is assigned
(that is, the outage was identified correctly) while dark red colors
indicate that the probability assigned to the true outage scenario is
less than $0.5$.  For example, the second bar from the left in
Figure~\subref*{fig:PMU.placement.GroupLASSO.3}, which corresponds to
two PMUs, corresponds to the following distribution of probabilities
assigned to the correct outage scenario, among the 3350 test
instances.
\begin{center}\footnotesize
\setlength{\tabcolsep}{.5em}
\begin{tabular}{|c||c|c|c|c|c|c|}\hline
Probability  & $[.9,1]$ & $[.8,.9]$ & $[.7,.8]$ & $[.6,.7]$ & $[.5,.6]$ & $[0,.5]$ \\\hline
\# of Instances &    963 &   144 &   135 &   161 &   206 &   1741  \\\hline
\end{tabular}
\end{center}
Note that the dark blue color occupies a fraction $963/3350$ of the
bar, medium blue occupies $144/3350$, and so on.

When only one PMU is installed, that bus naturally serves as the angle
reference, so no phase angle difference information is available, and
identification cannot be performed. As expected, identification becomes more
reliable as PMUs are installed on more buses. The value $\tau=.1$ (Figure
\subref*{fig:PMU.placement.GroupLASSO.3}) appears to select locations better
than the smaller choices of regularization parameter. For this value, about 10
buses are sufficient to assign a probability of greater than 90\% to the
correct outage event for more than 90\% of the test cases, while near-perfect
identification occurs when 30 PMUs are installed. Note that for $\tau=.1$,
there is only slow marginal improvement after 10 buses; we see a similar
pattern for the other values of $\tau$.  The locations added after the initial
selection are being chosen on the basis of information from the single
regularized problem \eqref{eq:MLE.group.l1}, so locations added later may be
providing only redundant information over locations selected earlier.

Figure~\ref{fig:PMU.placement.greedy} shows performance of the Greedy
Heuristic, plotted in the same fashion as in
Figure~\ref{fig:PMU.placement.GroupLASSO}.  For each value of $\tau$,
Algorithm~\ref{alg:PMU.placement.greedy} is performed with
$\cR^0=\emptyset$, with iterations continuing until there is no group
$s\in \cS \setminus \cR^{l-1}$ such that $q_s(\beta^l)>0$. Termination
occurs at 24, 16, and 11 PMU locations for the values $\tau=10^{-5},
10^{-3}$, and $10^{-1}$, respectively. As the value of $\tau$
increases, the number of PMUs which are selected for line outage
identification decreases. We can see by comparing
Figures~\ref{fig:PMU.placement.GroupLASSO} and
\ref{fig:PMU.placement.greedy} that classification performance
improves more rapidly as new locations are added for the Greedy
Heuristic than for the GroupLASSO Heuristic. Larger values of $\tau$
give slightly better results. We note
(Figure~\subref*{fig:PMU.placement.greedy.3} that almost perfect
identification occurs with only 16 PMU locations, while only 6
locations suffice to identify 90\% of outage events with high
confidence.

Although we can manipulate the GroupLASSO technique to achieve
sparsity equivalent to the Greedy Heuristic (by choosing a larger
value of $\tau$), the PMUs selected by the latter give much better
identification performance on this test set.  In
Table~\ref{tbl:Compare.GL.GR}, the parameter $\tau$ in the GroupLASSO
heuristic is chosen manually, to find the solutions with 10 PMUs and
15 PMUs for the 57 Bus system. Performance is compared to that
obtained from the Greedy Heuristic, with a much smaller value of
$\tau$. Results for the Greedy Heuristic are clearly superior.

\begin{table*}\centering
\caption{Comparison between GroupLASSO and Greedy Heuristic Selections on 57-Bus System}
\label{tbl:Compare.GL.GR}
\vspace{-5pt}
\begin{tabular}{|c|c|c|c||c|c|c||c|c|c|}\hline
\# of & \mr{2}{*}{Strategy} & \mr{2}{*}{$\tau$}
    & \mr{2}{*}{PMU Locations}
    & \mc{3}{c||}{Probability} & \mc{3}{c|}{Ranking} \\
\cline{5-10}
PMUs & & & & $\ge0.9$ & $\ge0.7$ & $\ge0.5$ & $1$ & $\le2$ & $\le3$\\\hline\hline
\mr{2}{*}{10} & GroupLASSO & $1.1$
                   & $1^*$ 8 17 27 28 51 52 53 54 55
                   & 72.8\% & 73.1\% & 78.8\% & 78.9\% & 92.7\% & 95.7\% \\
                   \cline{2-10}
                   & Greedy & $1.2\times10^{-1}$
                   & $1^*$ 2 17 19 26 39 40 45 46 57
                   & 92.6\% & 92.7\% & 94.3\% & 94.3\% & 99.7\% & 99.9\% \\
                   \hline\hline
\mr{2}{*}{15} & GroupLASSO & $8.0\times10^{-1}$
                   & $1^*$ 2 4 17 23 27 28 43 46 47 51 52 53 54 55
                   & 82.8\% & 82.8\% & 88.3\% & 88.3\% & 95.7\% & 95.8\% \\
                   \cline{2-10}
                   & Greedy & $1.7\times10^{-3}$
                   & $1^*$ 2 5 12 17 20 21 26 39 40 43 45 46 54 57
                   & 98.3\% & 98.3\% & 98.3\% & 98.3\% &  100\% &  100\% \\
                   \hline
\mc{10}{r}{$^*$ indicates the reference bus.}
\end{tabular}
\vspace*{-7pt}
\end{table*}

\subsubsection{Greedy Heuristic on 14, 30, 57 and 118 Bus System}

\begin{table*}\centering
\caption{Line Outage Detection Test Set with PMUs on About $\sim25\%$
  of Buses. \label{tbl:PMU.placement}}
\vspace{-5pt}
\begin{tabular}{|c|c|c|c||c|c|c||c|c|c|}\hline
\mr{2}{*}{System}  & \mr{2}{*}{$\tau$} & \# of
    & \mr{2}{*}{PMU Locations}
    & \mc{3}{c||}{Probability} & \mc{3}{c|}{Ranking} \\
\cline{5-10}
 & & PMUs & & $\ge0.9$ & $\ge0.7$ & $\ge0.5$ & $1$ & $\le2$ & $\le3$\\\hline\hline
\mr{2}{*}{ 14-Bus} & $5\times10^{-2}$ & 3
                   & $1^*$ 7 12
                   & 99.6\% & 99.7\% & 99.8\% & 99.8\% &  100\% &  100\% \\
                   \cline{2-10}
                   & $5\times10^{-3}$ & 3
                   & $1^*$ 11 12
                   &  100\% &  100\% &  100\% &  100\% &  100\% &  100\% \\
                   \hline\hline
\mr{2}{*}{ 30-Bus} & $5\times10^{-2}$ & 4
                   & $1^*$ 3 23 30
                   & 99.6\% & 99.6\% & 99.6\% & 99.6\% &  100\% &  100\% \\
                   \cline{2-10}
                   & $5\times10^{-3}$ & 5
                   & $1^*$ 3 14 22 29
                   &  100\% &  100\% &  100\% &  100\% &  100\% &  100\% \\
                   \hline\hline
\mr{2}{*}{ 57-Bus} & $5\times10^{-2}$ & 12
                   & $1^*$ 2 5 17 21 26 39 40 45 46 54 57
                   & 97.1\% & 97.1\% & 97.1\% & 97.1\% & 99.8\% & 99.8\% \\
                   \cline{2-10}
                   & $5\times10^{-3}$ & 14
                   & $1^*$ 2 5 17 20 21 26 39 40 41 45 46 54 57
                   & 98.5\% & 98.5\% & 98.5\% & 98.5\% & 99.9\% & 99.9\% \\
                   \hline\hline
\mr{3}{*}{118-Bus} & $5\times10^{-2}$ & 15
                   & 2 22 29 36 48 58 62 63 $69^*$ 81 91 95 106 108 115
                   & 94.2\% & 94.2\% & 94.2\% & 94.2\% & 96.2\% &  96.3\% \\
                   \cline{2-10}
                   & $5\times10^{-3}$ & 21
                   & \tcell{3 13 29 35 43 47 55 58 62 63 $69^*$\\ 75 81 82 91 93 104 106 107 113 115 119}
                   & 99.3\% & 99.5\% & 99.6\% & 99.6\% & 99.9\% &  99.9\% \\
                   \hline
\mc{10}{r}{$^*$ indicates the reference bus.}
\end{tabular}
\vspace*{-7pt}
\end{table*}

We applied the Greedy Heuristic to 14, 30, 57 and 118 Bus Systems with
two values of $\tau=5\times10^{-2}$ and $\tau=5\times10^{-3}$), and
found that the phasor measurements from the small set of buses are
enough to have the similar line outage identification performance to
the full measurement cases.  Table~\ref{tbl:PMU.placement} shows the
PMU locations selected for each case, and line outage identification
performance.  Identification performance is hardly degraded from the
fully instrumented case, even when phasor measurements are available
from only about $~25\%$ of buses.

\section{Conclusions}\label{sec:conclusion}

A novel approach to identify single line outage using MLR model is
proposed in this paper. The model employs historical load demand data
to train a multiclass logistic regression classifier, then uses the
classifier to identify outages in real time from streaming PMU data.
Numerical results obtained on IEEE 14, 30, 57 and 118 bus systems
prove that the approach can identify outages reliably.

With this line outage identification framework, the optimal placement
of PMU devices to identify the line outage is also
discussed. Heuristics are proposed to decide which buses should be
instrumented with PMUs. Experimental results show that detection is
almost as good  when just 25\% of buses are instrumented with PMUs
as when PMUs are attached to all buses.


\section*{Acknowledgment}
The authors thank Professor Chris DeMarco and Mr. Jong-Min Lim for
allowing the use of their synthetic 24-hour electric power demand data
sets in our experiments, and for valuable discussions and guidance on
this project.

\bibliographystyle{IEEEtran}
\bibliography{line_outage_identification}

\clearpage
\appendices
\pagenumbering{roman}

\section{Visualizing the Solution of the PMU Placement Problems}
\begin{figure}[h]\centering
\subfloat[IEEE 30 Bus System ($\tau=5\times10^{-3}$, 5 PMUs)]{
  \includegraphics[width=.45\textwidth]{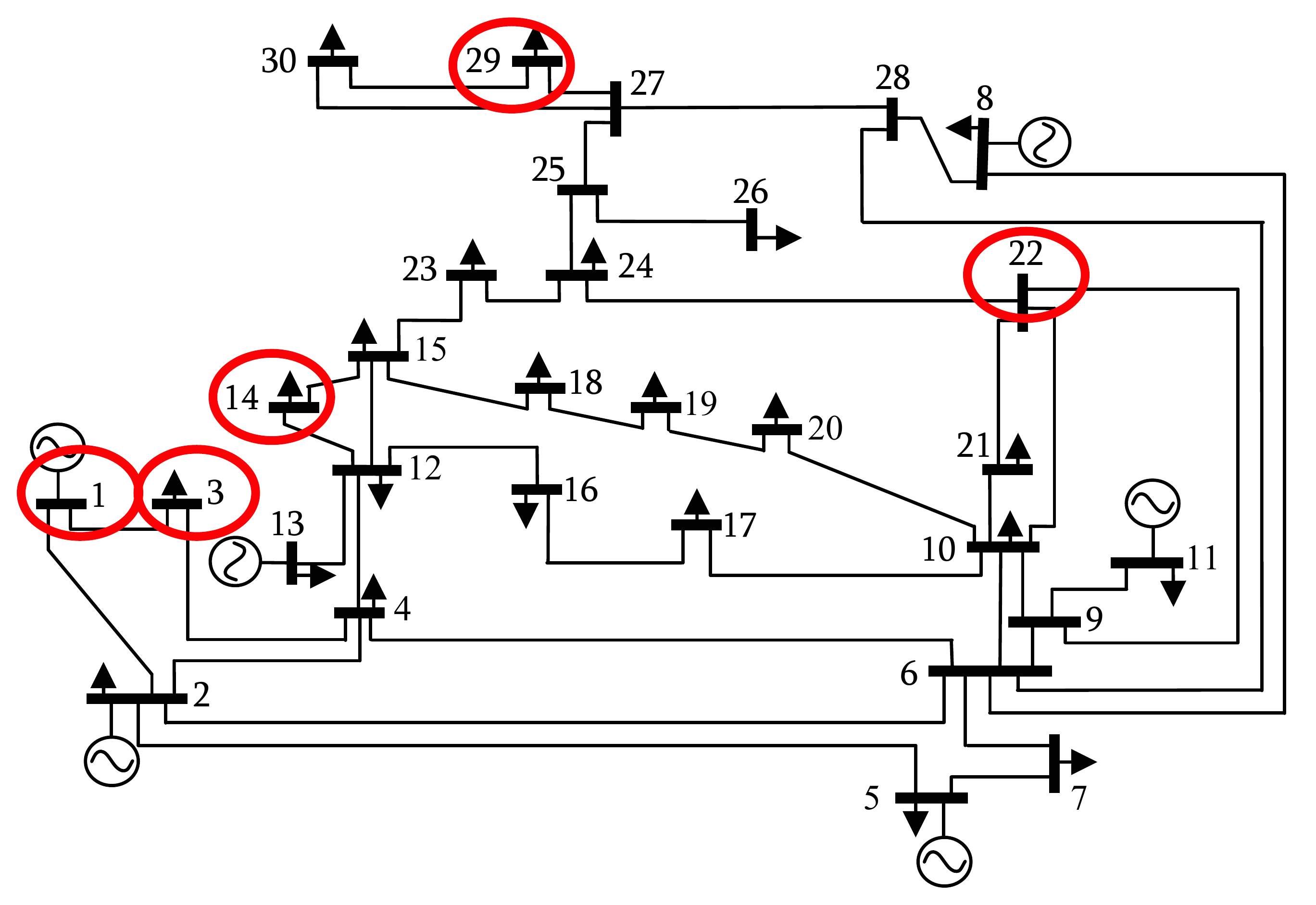}}\\
\subfloat[IEEE 57 Bus System ($\tau=5\times10^{-3}$, 14 PMUs)]{
  \includegraphics[width=.45\textwidth]{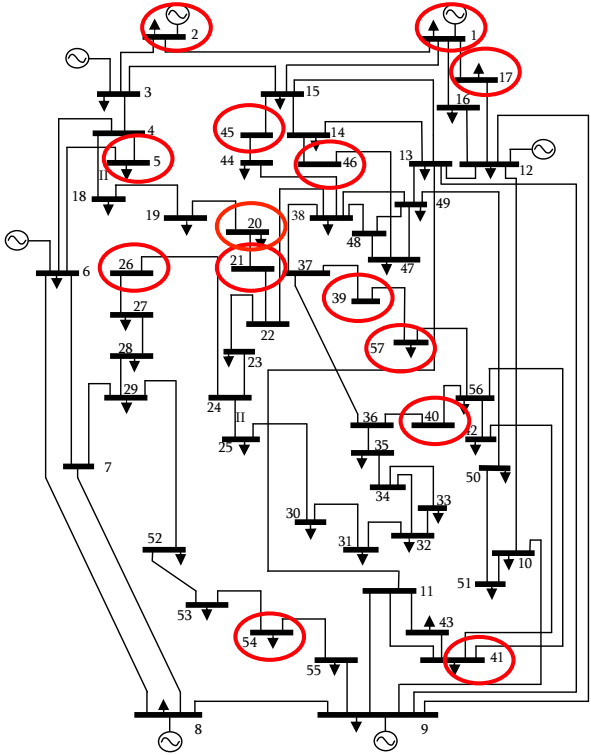}}
\caption{PMU Locations for IEEE 30 Bus and IEEE 57 Bus Systems.
  (System diagrams are taken from \cite{AllL08,GasA09})}\label{fig:PMU.location}
\end{figure}

The location of PMUs for the IEEE 30-Bus and IEEE 57-Bus Systems (from
Table~\ref{tbl:PMU.placement}) are displayed in Figure \ref{fig:PMU.location},
with instrumented buses indicated by red circles.

\section{Extension: Use of Explicit Line Outage
  Information}\label{sec:explicit.info}

We have assumed so far that only voltage angle and magnitude data from
PMUs is used in detecting line outages. In fact, PMUs provide other
information that is highly relevant for this purpose. For example,
when the PMU measures current of a particular line (incident on a bus) it can
detect immediately when an outage occurs on that line; we do not need
to rely in the indirect evidence of voltage changes at the other PMUs.
Another factor to consider is that when a decision is made to install
a PMU at a particular bus, it is conventional to measure {\em
  all lines} that are incident on that bus, as the marginal cost of doing so is
minimal. Although the phasor measurements are the same at all PMUs
near a single bus, each of these PMUs provides direct information
about the lines to which they are attached. Thus, if we choose to equip a
particular bus with PMUs, we can immediately detect outages on all
lines that touch that bus. In particular, if we install PMUs on {\em
  all} buses in the system, we have direct monitoring of all lines,
and the outage detection problem becomes trivial.

We can extend the multiclass logistic regression technique to make
effective use of these direct observations in choosing optimal buses
for PMU placement. The key modification is to extend each feature
vector $\beta_k$ to include additional entries that indicate the buses
that touch line $k$. The observation vectors and the groups $\cP_s$,
$s=1,2,\dotsc,N$ are extended correspondingly.

For each line $k=1,2,\dotsc,K$, let us define the following
quantities:
\[
  \cT_k := \set{t_k^1, t_k^2}  \subset \set{1,2,\dotsc,N }
\]
where $t_k^1$ and $t_k^2$ are indices of buses touched by line $k$.

We extend each observation vector $\ub{X}$ by appending $2K$
additional elements to form $\ub{\ub{X}}$, where each such vector has
the form
\beq\label{eq:Xbarbar}
\ub{\ub{X}}:=\bm{\ub{X} \\ {L_k}}
\eeq
for some $k=1,2,\dotsc,K$, where $L_k$ is the $k$th column of the
$2K\times K$ matrix $L$ defined as follows, for some $\eta>0$:
\beq \label{def:L}
 L := \bm{ \eta & 0 & \cdots & 0 \\
           \eta & 0 & \cdots & 0 \\
           0 & \eta & \cdots & 0 \\
           0 & \eta & \cdots & 0 \\
           \vdots & \vdots & \ddots & \vdots \\
           0 & 0 & \cdots & \eta \\
           0 & 0 & \cdots & \eta }.
\eeq
The two nonzero entries in each column $L_k$ indicate which two buses
can detect outage of line $k$ directly.  In other words, if line $k$
fails, we flag the PMUs on the buses that touch that line with a value
$\eta$, since a fault on line $k$ is immediately detectable from the
buses $t^1_k$ and $t^2_k$. We can say that the first part of the
combined observation vector $\ub{\ub{X}}$ contains {\em indirect}
(voltage phasor, $\ub{X}$) observations while the second part contains
{\em direct} (line outage, $L_k$) observations.

We need to extend too the definition \eqref{eq:Ps} of the groups $\cP_s$,
$s=1,2,\dotsc,N$.  We now distribute the additional $2K$ entries in the feature
vector to these groups. The additional entries associated with bus $s$ are
those in the index set $\cD_s$ defined as follows:
\begin{align*}
 \cD_s & = \setc{2(k-1)+i}{t_k^i=s \;\; \mbox{for $k=1,2,\cdots,K$, $i=1,2$}}.
\end{align*}
For the combined observation vector $\ub{\ub{X}}$, we define groups
$\cP_s'$, $s=1,2,\cdots,N$ :
\[
  \cP_s'     = \cP_s \cup \setc{2N+2+d}{d\in\cD_s}\;\;\mbox{for $s=1,2,\cdots,N$},
\]

If PMUs are installed on every line, direct observations will identify
each outage perfectly, so the solution of the maximum likelihood
problem is rather trivial. In approximate solutions to the problem,
the weight vector $\beta_k$ for line $k$ will have large positive
entries in positions $2N+2+i$ for which $L_{ik}=\eta$, and large
negative entries in positions $2N+2+i$ for which $L_{ik}=0$.  This
would yield $\beta_k^T\ub{\ub{X}}$ large and positive for observation
vectors $\ub{\ub{X}}$ that indicate a line-$k$ outage, with
$\beta_k^T\ub{\ub{X}}$ large and negative if there is no outage on
line $k$, leading to assigned probabilities close to $1$ and $0$,
respectively. (Entries in $\beta_k$ corresponding to the indirect
observations may also have meaningfully large values, but these are
less significant in the completely observed case.)

When PMUs are installed on a subset of buses, outages on some lines
will be observed only {\em indirectly}, so the indirect observations
in components $i=1,2,\dotsc,2N+2$ of the vector $\ub{\ub{X}}$ are
critical to identification performance on those lines that are not
directly observed.

We incorporate direct observations into our outage identification
strategy in the following ways.
\begin{itemize}
  \item Indirect. Direct observations are ignored. We use only the
    observation vector $\ub{X}$, as in Section~\ref{sec:result}.
  \item Combined (Direct+Indirect). Direct observations are
    incorporated into the observation vector, and we so MLR
    classification with the vectors $\ub{\ub{X}}$.
  \item Prescreening. Instead of including the direct observation in
    the observation vector, line outages that can be identified by the
    direct observation are screened out {\em before} the MLR is
    applied. The number of outcomes in MLR is reduced since we do not
    need to consider the line outages identified already by the direct
    observation. Observation vectors $\ub{X}$ are used to train the
    MLR for those line outages that are not observed directly.
  \item Postscreening. First, we train an MLR classifier using only the
    indirect observations in vector $\ub{X}$. Then, during
    classification, we override the prediction result from the MLR
    when a direct observation is available for the line in question.
    Note that results from this strategy cannot be worse than results
    for the Indirect strategy.
\end{itemize}

We also compare solutions of the PMU placement problem using the
Indirect strategy (as in Section~\ref{sec:result}) and the Combined
strategy. We solve these two variants of the placement problem for the
57-bus case with the Greedy Heuristic of
Algorithm~\ref{alg:PMU.placement.greedy}, setting $\tau=10^{-2}$ and
the number of PMUs $r$ to the values $5$ and $10$.
We note that the reference bus is always selected as one of the PMU
locations, and it is used only to provide the phase angle reference
for all the strategies above, as in the Indirect case. (In using the
reference bus PMU in this restricted way, we allow a fairer comparison
between the Indirect strategy and the strategies that use direct
observations.)

\begin{table*}\centering
\caption{Use of Explicit Line Outage Information
  ($\tau=10^{-2}$, $^*$ indicates the reference bus.)}
\label{tbl:ext.DO}
\subfloat[PMU Placement Based on Indirect Observations]{
  \label{tbl:ext.DO.ind}
  \shortstack[r]{
  \begin{tabular}{|c||c|c|c||c|c|c||c|c|c||c|c|c|}\hline
  \mr{3}{*}{Strategy}
    & \mc{6}{c||}{$1^*$ 5 20 21 57 (5 PMUs)}
    & \mc{6}{c|}{$1^*$ 5 20 21 26 39 40 43 54 57 (10 PMUs)}\\
  \cline{2-13}
    & \mc{3}{c||}{Probability} & \mc{3}{c||}{Ranking}
    & \mc{3}{c||}{Probability} & \mc{3}{c|}{Ranking}\\
  \cline{2-13}
   & $\ge0.9$ & $\ge0.7$ & $\ge0.5$ & $1$ & $\le2$ & $\le3$
   & $\ge0.9$ & $\ge0.7$ & $\ge0.5$ & $1$ & $\le2$ & $\le3$\\\hline\hline
  Indirect   & 83.1\% & 83.2\% & 88.0\% & 88.0\% & 99.4\% & 99.8\%
             & 94.1\% & 94.1\% & 94.1\% & 94.1\% & 99.1\% & 99.1\% \\
  Combined   & 83.0\% & 83.8\% & 86.4\% & 86.4\% & 97.5\% & 99.4\%
             & 95.4\% & 95.4\% & 95.4\% & 95.4\% & 98.8\% & 98.8\% \\
  Prescreening & 84.8\% & 85.1\% & 88.1\% & 88.1\% & 98.4\% & 99.6\%
                & 95.2\% & 95.2\% & 95.3\% & 95.3\% & 99.8\% & 99.9\% \\
  Postscreening & 86.1\% & 86.2\% & 89.5\% & 89.5\% & 99.4\% & 99.8\%
                 & 94.2\% & 94.2\% & 94.2\% & 94.2\% & 99.1\% & 99.1\% \\
  \hline
  \end{tabular}\\
  \begin{tabular}{|c|ccccccccc|c|}\hline
  Selected Bus     &  5 & 20 & 21 & 26 & 39 & 40 & 43 & 54 & 57 & Total \\\hline
  \# of Lines Touching the Bus &  2 &  2 &  2 &  2 &  2 &  2 &  2 &  2 &  2 & 18 \\\hline
  \end{tabular}}
}\\
\subfloat[PMU Placement Based on Combined (Direct + Indirect) Observations ($\eta=1$)]{
  \label{tbl:ext.DO.dir.1}
  \shortstack[r]{
  \begin{tabular}{|c||c|c|c||c|c|c||c|c|c||c|c|c|}\hline
  \mr{3}{*}{Strategy}
    & \mc{6}{c||}{$1^*$ 6 9 12 56 (5 PMUs)}
    & \mc{6}{c|}{$1^*$ 6 9 12 15 22 39 49 54 56 (10 PMUs)}\\
  \cline{2-13}
    & \mc{3}{c||}{Probability} & \mc{3}{c||}{Ranking}
    & \mc{3}{c||}{Probability} & \mc{3}{c|}{Ranking}\\
  \cline{2-13}
   & $\ge0.9$ & $\ge0.7$ & $\ge0.5$ & $1$ & $\le2$ & $\le3$
   & $\ge0.9$ & $\ge0.7$ & $\ge0.5$ & $1$ & $\le2$ & $\le3$\\\hline\hline
  Indirect   & 77.2\% & 77.4\% & 83.1\% & 83.2\% & 97.4\% & 99.0\%
             & 83.7\% & 83.7\% & 87.5\% & 87.5\% & 93.5\% & 93.5\% \\
  Combined   & 83.1\% & 83.1\% & 87.4\% & 87.4\% & 92.7\% & 92.8\%
             & 93.1\% & 93.1\% & 94.7\% & 94.7\% & 98.4\% & 98.5\% \\
  Prescreening & 83.3\% & 83.8\% & 88.2\% & 88.2\% & 93.4\% & 93.4\%
                & 93.8\% & 93.8\% & 95.2\% & 95.2\% & 98.4\% & 98.4\% \\
  Postscreening & 82.7\% & 82.8\% & 87.2\% & 87.2\% & 98.5\% & 99.6\%
                 & 93.7\% & 93.7\% & 95.1\% & 95.1\% & 97.4\% & 97.4\% \\
  \hline
  \end{tabular}\\
  \begin{tabular}{|c|ccccccccc|c|}\hline
  Selected  Bus     &  6 &  9 & 12 & 15 & 22 & 39 & 49 & 54 & 56 & Total \\\hline
  \# of Lines Touching the Bus &  4 &  6 &  5 &  5 &  3 &  2 &  4 &  2 &  4 & 35 \\\hline
  \end{tabular}}
}\\
\subfloat[PMU Placement Based on Combined (Direct + Indirect) Observations ($\eta=10^{-2}$)]{
  \label{tbl:ext.DO.dir.2}
  \shortstack[r]{
  \begin{tabular}{|c||c|c|c||c|c|c||c|c|c||c|c|c|}\hline
  \mr{3}{*}{Strategy}
    & \mc{6}{c||}{$1^*$ 5 9 49 56 (5 PMUs)}
    & \mc{6}{c|}{$1^*$ 5 9 21 26 39 45 46 49 56 (10 PMUs)}\\
  \cline{2-13}
    & \mc{3}{c||}{Probability} & \mc{3}{c||}{Ranking}
    & \mc{3}{c||}{Probability} & \mc{3}{c|}{Ranking}\\
  \cline{2-13}
   & $\ge0.9$ & $\ge0.7$ & $\ge0.5$ & $1$ & $\le2$ & $\le3$
   & $\ge0.9$ & $\ge0.7$ & $\ge0.5$ & $1$ & $\le2$ & $\le3$\\\hline\hline
  Indirect   & 79.0\% & 79.3\% & 84.3\% & 84.5\% & 98.9\% & 99.8\%
             & 93.3\% & 93.3\% & 94.7\% & 94.7\% & 99.4\% & 99.9\% \\
  Combined   & 85.4\% & 85.4\% & 90.0\% & 90.0\% & 95.5\% & 95.6\%
             & 97.4\% & 97.4\% & 97.4\% & 97.4\% & 99.9\% & 99.9\% \\
  Pre-Screening & 84.5\% & 84.5\% & 88.5\% & 88.5\% & 98.4\% & 98.9\%
                & 96.1\% & 96.1\% & 96.1\% & 96.1\% & 97.9\% & 97.9\% \\
  Post-Screening & 85.1\% & 85.4\% & 89.5\% & 89.7\% & 99.7\% &  100\%
                 & 98.0\% & 98.0\% & 98.0\% & 98.0\% &  100\% &  100\% \\
  \hline
  \end{tabular}\\
  \begin{tabular}{|c|ccccccccc|c|}\hline
  Selected Bus     &  5 &  9 & 21 & 26 & 39 & 45 & 46 & 49 & 56 & Total \\\hline
  \# of Lines Touching the Bus &  2 &  6 &  2 &  2 &  2 &  2 &  2 &  4 &  4 & 26 \\\hline
  \end{tabular}}
}
\vspace{-3pt}
\end{table*}

Experimental results using PMU placements based on Indirect and Combined
observations, and using each of the four classification strategies described
above, are shown in Table~\ref{tbl:ext.DO}, using a similar format to
Tables~\ref{tbl:Compare.GL.GR} and \ref{tbl:PMU.placement}.  When the PMU
locations are selected using only indirect observations
(Table~\subref*{tbl:ext.DO.ind}), the advantage of using direct line outage
information during classification is not significant, especially when the
larger number of $10$ PMUs is installed.  This observation is not too
surprising.  The biggest voltage phasor changes are produced by outages that
are close to a bus, so even when a line outage is not detected by direct
observation, it can usually be detected reliably by its ``indirect'' effect on
nearby buses.

In Tables~\subref*{tbl:ext.DO.dir.1} and \subref*{tbl:ext.DO.dir.2},
the PMU locations are selected on the basis of the combined
vectors. When only indirect data is used during classification,
results are much worse, as the locations have been chosen under the
assumption that direct observation data will be available. In fact,
the results in Table~\subref*{tbl:ext.DO.dir.1} are generally slightly
worse than those of Table~\subref*{tbl:ext.DO.ind}. This is again
because too much reliance is placed on direct observation in selecting
PMU locations, and detection power is diminished slightly for those
outages that are detectable only indirectly. Note that the PMU
locations in Table~\subref*{tbl:ext.DO.dir.1} are essentially those
with the greatest numbers of lines connected: a total of 35 in
Table~\subref*{tbl:ext.DO.dir.1} (for $r=10$), as compared with 18 in
Table~\subref*{tbl:ext.DO.ind}.

To reduce the weight placed on direct information in PMU placement, we
scale down the values $\eta$ in \eqref{def:L}. Reducing $\eta$ from
$1$ to $10^{-2}$ appears to strike a better balance between the use of
direct and indirect information. Table~\subref*{tbl:ext.DO.dir.2}
shows a marked improvement over Table~\subref*{tbl:ext.DO.ind} (which
weights the direct observations more heavily), and slight improvements
by most measures over Table~\subref*{tbl:ext.DO.dir.1}, which uses
only indirect information. The total number of lines that are directly
connected to buses with PMUs (and which can thus be observed directly)
in Table~\subref*{tbl:ext.DO.dir.1} is about halfway between the
corresponding statistics in Tables~\subref*{tbl:ext.DO.ind} and
\subref*{tbl:ext.DO.dir.2}.

\section{Implementation of \SpaRSA{}}

We solve the regularized convex optimization problem \eqref{eq:MLE.p}
with the \SpaRSA{} algorithm \cite{WriNF08}, a simple first-order
approach that exploits the structure. We briefly describe the approach
here, referring to \cite{WriNF08} for further details.

The \SpaRSA{} subproblem of \eqref{eq:MLE.p} at iteration $n$ is
defined as follows, for some scalar parameter $\alpha_n\in
\mathbb{R}^+$:
\beq\label{eq:SpaRSA.sub}
  \beta^{n+1} := \oargmax{\beta}
      \frac{1}{2}\norm{\beta-\gamma^n}_{F}^2
     + \tau\frac{1}{\alpha_n}w(\beta)
\eeq
where $\|\cdot\|_F$ is the Frobenius norm of a matrix, $\gamma^n :=
\bm{\gamma_1^n & \gamma_2^n & \dotsc & \gamma_K^n}$ with
\[
  \gamma_k^n := \beta_k^n - \frac{1}{\alpha_n}\nabla_{\beta_k}f(\beta),
\]
and
\[
  \nabla_{\beta_i}f(\beta) =
    \sum_{p:i=y_p}x_p - \sum_{j=1}^M\frac{x_je^{\ip{\beta_i}{x_j}}}
      {\sum_{k=1}^Ke^{\ip{\beta_k}{x_j}}}.
\]
When no regularization term is present ($w(\beta)=0$), the solution
for \eqref{eq:SpaRSA.sub} is $\beta^{n+1}=\gamma^n$, so the approach
reduces to the steepest descent algorithm on $f$ with step length
$1/\alpha_n$.

In the PMU placement problem, our regularizer $w_\cS(\beta)$ is
group-separable. Thus the subproblem \eqref{eq:SpaRSA.sub} can be
divided into independent problems of the form
\[
  [\beta^{n+1}]_{\cP_s} := \oargmax{\hat{\beta}}
     \frac{1}{2}\norm{[\beta]_{\cP_s}-[\gamma^n]_{\cP_s}}_2^2
     + \tau\frac{1}{\alpha_n}q_s(\beta),
\]
for all $s \in \cS$, where (as defined above), $[A]_{\cP_s}$ is the
submatrix of $A$ consisting of the rows whose indices are in
${\cP_s}$.  Since the penalty function $q_s(\beta)$ is the
$\ell_2$-norm, this subproblem has a closed form solution
\cite{WriNF08}, as follows:
\[
  [\beta^{n+1}]_{\cP_s} = [\gamma^n]_{\cP_s}
       \frac{\max\set{\norm{[\gamma^n]_{\cP_s}}_2-\tau\alpha_n^{-1},0}}
            {\max\set{\norm{[\gamma^n]_{\cP_s}}_2-\tau\alpha_n^{-1},0}+\tau\alpha_n^{-1}}.
\]
For any row $i$ of $\beta$ that does not belong to any $\cP_s$, we
have simply $[\beta^{n+1}]_i=[\gamma^n]_i$.  Different strategies can
be used to choose $\alpha_n$.  We increase $\alpha_n$ at each
iteration until sufficient decrease is obtained in the objective,
terminating when $\alpha_n$ grows too large (indicating that a
solution is nearby).

\end{document}